\documentclass[12pt]{article}

\usepackage{multirow, multicol}
\usepackage{array} 
\usepackage{paralist} 
\usepackage{verbatim} 

\usepackage{amsmath, amssymb, enumerate}
\usepackage{pictexwd,dcpic}
\usepackage[final]{graphicx}

\newtheorem{theorem}{Theorem}
\newtheorem{corollary}[theorem]{Corollary}
\newtheorem{remark}[theorem]{Remark}
\newtheorem{lemma}[theorem]{Lemma}
\newtheorem{definition}[theorem]{Definition}
\newtheorem{proposition}[theorem]{Proposition}
\newtheorem{example}[theorem]{Example}

\newcommand{\pend}{\hfill $\square$}

\numberwithin{equation}{section}  
\numberwithin{figure}{section}    
\numberwithin{theorem}{section}

\newcommand{\of}[1]{\ensuremath{\left( #1 \right)}}
\newcommand{\norm}[1]{\ensuremath{\left\| #1 \right\|}}
\newcommand{\abs}[1]{\ensuremath{\left| #1 \right|}}
\newcommand{\cb}[1]{\ensuremath{ \left\{ #1 \right\} }}
\newcommand{\sqb}[1]{\ensuremath{ \left[ #1 \right] }}

\newcommand{\bs}{\backslash}

\newcommand{\eps}{\ensuremath{\varepsilon}}

\newcommand{\vp}{\ensuremath{\varphi}}

\newcommand{\R}{\mathrm{I\negthinspace R}}
\newcommand{\OLR}{\overline{\mathrm{I\negthinspace R}}}
\newcommand{\N}{\mathrm{I\negthinspace N}}
\newcommand{\E}{{\mathbb E}}

\renewcommand{\P}{\ensuremath{\mathcal{P}}}

\newcommand{\G}{\ensuremath{\mathcal{G}}}
\newcommand{\F}{\ensuremath{\mathcal{F}}}

\newcommand{\A}{\ensuremath{\mathcal{A}}}

\newcommand{\Min}{{\rm Min\,}}

\newcommand{\WEff}{{\rm wEff}}

\newcommand{\Sol}{{\rm Sol\,}}

\newcommand{\dom}{{\rm dom \,}}

\newcommand{\cl}{{\rm cl \,}}

\newcommand{\co}{{\rm co \,}}

\newcommand{\cone}{{\rm cone\,}}

\newcommand{\Int}{{\rm int\,}}

\newcommand{\dc}{{\rm dcl\,}}

\newcommand{\lev}{{\rm lev }}

\newcommand{\lel}{\preccurlyeq}
\newcommand{\leu}{\curlyeqprec}

\newcommand{\triup}{{\rm \vartriangle}}
\newcommand{\trido}{{\rm \triangledown}}

\newcommand{\One}{\mathrm{1\negthickspace I}}

\usepackage[colorlinks=true, linkcolor=blue, citecolor=blue, urlcolor=blue, filecolor=blue]{hyperref}

\usepackage{xcolor}
\definecolor{color0}{gray}{.50}
\definecolor{color1}{rgb}{0,.2,.8}
\definecolor{color2}{rgb}{1,.2,0}
\definecolor{color3}{rgb}{.8,.5,1}
\newcommand{\f}{\color{color1}}

\begin{document}

\title{Set Relations and Approximate Solutions in Set Optimization}

\author{Giovanni Crespi\footnote{Universit{\`a} degli Studi dell'Insubria, Department of Economics, \href{mailto:giovanni.crespi@uninsubria.it}{giovanni.crespi@uninsubria.it}}, Andreas H. Hamel\footnote{Free University of Bozen-Bolzano, Faculty of Economics and Management, \href{mailto:andreas.hamel@unibz.it}{andreas.hamel@unibz.it}}, Matteo Rocca\footnote{Universit{\`a} degli Studi dell'Insubria, Department of Economics, \href{mailto:matteo.rocca@uninsubria.it}{matteo.rocca@uninsubria.it}}, Carola Schrage\footnote{Free University of Bozen-Bolzano, Faculty of Economics and Management, \href{mailto:carola.schrage@unibz.it}{carola.schrage@unibz.it}}}

\maketitle
\begin{abstract} Via a family of monotone scalar functions, a preorder on a set is extended to its power set and then used to construct a hull operator and a corresponing complete lattice of sets. A function mappping into the preordered set is extended to a complete lattice-valued one, and concepts for exact and approximate solutions for corresponding set optimization problems are introduced and existence results are given. Well-posedness for complete lattice-valued problems is introduced and characterized. The new approach is compared to existing ones in vector and set optimization, and its relevance is shown by means of many examples from vector optimization, statistics and mathematical economics \& finance.
\end{abstract}

\medskip\noindent
{\bf Keywords} set relation, hull operator, complete lattice, set optimization, approximate solution, well-posedness 

\medskip\noindent
{\bf Mathematics Subject Classification} 49J53, 47N10, 46N10

\medskip\noindent
{\bf Acknowledgement.} The work of all four coauthors was part of the project ``An abstract convexity approach to scalarization in set/vector optimization and to multi-utility maximization" funded by Free University of Bozen-Bolzano, Italy.

\section{Introduction}

A preordered set $(Z, \preceq)$ and a family of extended real-valued functions $\Psi$ on $Z$ which are monotone with respect to $\preceq$ are the two basic ingredients for the theory developed in this paper. Such a family can be understood as a collection of ``elementary functions" as in abstract convexity (e.g., \cite{Rubinov13Book}), as a ``multi-utility representation" of an incomplete preference as in econonmics (e.g., \cite{EvrenOk11JME}, \cite{Evren14JET}), as a family which defines (e.g., stochastic dominance orders) or characterizes (e.g., a vector order via the bipolar theorem) an order relation. Since there is very little structure on $(Z, \preceq)$, it is difficult to effectively deal with optimization problems for functions mapping into $Z$, even to define reasonable solution concepts.

We propose a new approach which extends the order $\preceq$ on $Z$ to one on its power set, provides a corresponding hull operator which in turn admits the construction of a complete lattice of sets which are closed with respect to this hull operator. The complete lattice structure then admits to define solution concepts for corresponding optimization problems which--in contrast to the vast majority of papers in vector or set optimization--also involve the infimum or supremum as meaningful concepts.

These constructions are given through the family of monotone functions which are not seen as ``scalarizations," but as a mean to turn an optimization problem for a function mapping into a preordered set into one for a function mapping into a complete lattice of sets. 

Our approach can be understood as turning the tables with respect to scalarization procedures in vector/set optimization: instead of selecting a particular scalarization as a substitute for a vector- or set-valued objective, a family of scalarization functions is used, which defines or characterizes the order relation, to construct a set order relation and a corresponding complete lattice of sets. Usually, a multitude of representations of a given order on $Z$ via a family of scalar functions exists, some of them more useful than others. This gives the task to the decision maker to carefully select the ``right" family (depending on her/his purpose and the features of the mathematical model) before the optimization procedure even starts. The resulting hull operators and lattices may also be different.

Major contributions of this paper include the construction of set relations and complete lattices of sets via families of monotone scalar functions, new concepts for exact and approximate solutions involving the infimum or supremum for complete lattice-valued functions, corresponding Weierstrass type theorems and well-posedness results. These concepts also have a computational aspect: the analysis of algorithms for vector/set optimization problems such as \cite{LoehneSchrage13Opt} should be based on a clear understanding of what is considered a (good approximate) solution.

Many examples from different applied fields such as multi-criteria decision making, statistics, mathematical finance and economics show the relevance of our concepts; some of them are very much under discussion in their corresponding communities such as the multi-utility maximization problem for which basically no theory exists, or quantiles and stochastic dominance orders for multivariate random variables; some others are known, but have never been investigated from a complete lattice point of view such as the first and second stochastic dominance orders for the univariate random variables.

By the way of conclusion, the fundamental proposal of this paper is to replace a vector- or set-valued optimization problem by its complete lattice-valued extension which is considered to be a ``true" set optimization problem; in doing so one can obtain many concepts and results in striking parallelity to scalar results which is not possible in a ``pure" vector optimization setting.

The next section provides the complete lattice framework, new set order relations and hull operators, Section \ref{SecSolConcepts} includes the definitions of (approximate) solutions, and existence theorems as well as well-posedness results are given in Section \ref{SecExWellPo}. Throughout the paper, many examples are discussed emphazising features of the new concepts and comparing or linking them with existing ones. Quite a few of those examples are taken from the recent literatur in economics, finance, statistics and vector or set optimization.


\section{Set relations generated by families of scalar functions}

\label{SecSetRelations}

Let $(Z, \preceq)$ be a preordered set, i.e., $\preceq$ is a reflexive and transitive relation on $Z$. A function $\psi \colon Z \to \OLR:= \R\cup\cb{\pm\infty}$ is called monotone (with respect to $\preceq$) if
\[
z_1 \preceq z_2 \quad \Rightarrow \quad \psi(z_1) \leq \psi(z_2).
\]
Let $\P(Z)$ denote the set of all subsets of $Z$ (including $\emptyset$), i.e., the power set of $Z$. The function $\psi^\triup \colon \P(Z) \to \OLR$ defined by
\begin{equation}
\label{EqPsiInfExt}
\psi^\triup(D) = \inf_{z \in D}\psi(z)
\end{equation}
is called the {\em inf-extension} of $\psi \colon Z \to \OLR$ from $Z$ to $\P(Z)$ where $\inf_{z \in \emptyset}\psi(z) = +\infty$ by definition.

Let $\Psi$ be a family of monotone functions $\psi \colon Z \to \OLR$ and define a relation $\preceq_\Psi$ on $\P(Z)$ by
\[
D_1 \preceq_\Psi D_2  \quad \Leftrightarrow \quad \forall \psi \in \Psi \colon \psi^\triup(D_1) \leq \psi^\triup(D_2).
\]
The monotonicity of the functions $\psi$ guarantees 
\begin{equation}
\label{EqExtOrder}
z_1 \preceq z_2 \quad \Rightarrow \quad \cb{z_1} \preceq_\Psi \cb{z_2}.
\end{equation}
In this sense, $\preceq_\Psi$ is an extension of $\preceq$  from $Z$ to $\P(Z)$. Apparently, $\preceq_\Psi$ is a reflexive and transitive relation on $\P(Z)$. This procedure can be considered as a new way of introducing ``set relations,'' see \cite{KuroiwaTanakaTruong97NA, JahnTruong11JOTA, Hamel05Habil} and the survey \cite{HamelEtAl15Incoll} for more details and references.

Moreover, \eqref{EqExtOrder} also shows that the restriction of $\preceq_\Psi$ to $Z$ is an extension of $\preceq$ on $Z$. If the opposite implication in \eqref{EqExtOrder} is also true, i.e.,
\begin{equation}
\label{EqRepresenting}
z_1 \preceq z_2  \quad \Leftrightarrow \quad \cb{z_1} \preceq_\Psi \cb{z_2}
\quad \Leftrightarrow \quad \forall \psi \in \Psi \colon \psi(z_1) \leq \psi(z_2)
\end{equation}
the family $\Psi$ is called {\em representing} for $\preceq$. Below, it is shown that every preorder has a representing family of extended real-valued functions, but, unless explicitly stated, it is not assumed in the sequel that $\Psi$ is representing.

For each set $D \in \P(Z)$, we define
\[
\cl_\Psi D = \bigcap_{\psi \in \Psi}\cb{z \in Z \mid \psi^\triup(D) \leq \psi(z)}.
\]

The following lemma contains a few simple properties.

\begin{lemma} 
\label{LemExtToLattice}
Let $D, E \in \P(Z)$. Then

(a) $D \preceq_\Psi E$ if, and only if, $\cl_\Psi D \supseteq \cl_\Psi E$.

(b) $D \preceq_\Psi E$, $E \preceq_\Psi D$ if, and only if, $\cl_\Psi D = \cl_\Psi E$.

(c) it holds
\begin{equation}
\label{EqInfHull}
\psi^\triup(\cl_\Psi D) = \psi^\triup(D).
\end{equation}
\end{lemma}

{\sc Proof.} (a) Straightforward using the definitions of $\preceq_\Psi$ and $\cl_\Psi$. (b) is a consequence of (a).  (c) The inequality ``$\leq$" in \eqref{EqInfHull} follows from the obvious fact $D \subseteq \cl_\Psi D$. For the converse, observe $\psi(z) \geq \psi^\triup(D)$ for all $z \in \cl_\Psi D$ by definition of $\cl_\Psi D$, hence ``$\geq$." \pend

\begin{proposition}
\label{PropPsiHull}
The mapping $D \mapsto \cl_\Psi D$ is a hull operator, i.e., (i) $D \subseteq \cl_\Psi D$, (ii) $D \subseteq E$ implies $\cl_\Psi D \subseteq \cl_\Psi E$ and (iii) $\cl_\Psi D = \cl_\Psi(\cl_\Psi D)$ for all $D, E \in \P(Z)$. 
\end{proposition}

{\sc Proof.} (i) is straightforward from the definitions of $\psi^\triup$ and $\cl_\Psi D$. (ii) follows since $D \subseteq E$ implies $\psi^\triup(D) \geq \psi^\triup(E)$ for all $\psi \in \Psi$. (iii) is a consequence of \eqref{EqInfHull} with $D$ replaced by $\cl_\Psi D$. \pend

\medskip The previous proposition guarantees that the set
\begin{equation}
\label{EqPsiPower}
\P(Z, \Psi) = \cb{D \in \P(Z) \mid D = \cl_\Psi D}
\end{equation}
is well-defined. Indeed, we have $\cl_\Psi D \in \P(Z, \Psi)$ for all $D \in \P(Z)$ by (iii) of Proposition \ref{PropPsiHull}.

\begin{proposition} 
\label{PropPsiLattice} 
For $D, E \in \P(Z, \Psi)$, 
\[
D \preceq_\Psi E \quad \Leftrightarrow \quad D \supseteq  E.
\]
Moreover, the pair $\of{\P(Z, \Psi), \supseteq}$ is a complete lattice, and for each $\A \subseteq \P(Z, \Psi)$,
\[
\inf \A = \cl_\Psi \bigcup_{A \in \A} A \quad \text{and} \quad \sup \A = \bigcap_{A \in \A} A
\]
where the infimum and the supremum are taken with respect to $\supseteq$.
\end{proposition}

{\sc Proof.} The coincidence of $\preceq_\Psi$ and $\supseteq$ on $\P(Z, \Psi)$ is Lemma \ref{LemExtToLattice} (a). 

The formula for the supremum follows because the intersection of elements of $\P(Z, \Psi)$ again is an element of $\P(Z, \Psi)$ since $\cl_\Psi$ is a hull operator. Finally, we shall show
\begin{equation}
\label{EqInfSet}
\inf \A = \cl_\Psi \bigcup_{A \in \A} A = \bigcap_{\psi \in \Psi}\cb{z \in Z \mid \inf_{A \in \A}\psi^\triup(A) \leq \psi(z)}.
\end{equation}

Let, for the moment, denote $B = \cl_\Psi \bigcup_{A \in \A} A$. First, it needs to be shown that $B \supseteq A$ for all $A \in \mathcal A$. This follows from (i) of Proposition \ref{PropPsiHull}. 

Secondly, it must be verified that if $C \in \P(Z, \Psi)$ satisfies
\[
\forall A \in \mathcal A \colon C \supseteq A,
\]
then $C \supseteq B$ is true. Indeed, for such a $C$
\[
\forall \psi \in \Psi, \; \forall A \in \mathcal A \colon \psi^\triup(C) \leq \psi^\triup(A)
\]
follows, hence
\[
\forall \psi \in \Psi \colon \psi^\triup(C) \leq \inf_{A \in \mathcal A} \psi^\triup(A).
\]
If it can be shown that 
\[
\forall \psi \in \Psi \colon \psi^\triup(B) = \inf_{A \in \mathcal A} \psi^\triup(A),
\]
then $C \supseteq B$ follows from the definition of $\preceq_\Psi$ and (a). The above equality, which also produces the right equation in \eqref{EqInfSet}, is proved in the following lemma. \pend

\begin{lemma}
\label{LemInfScalar}
If $\mathcal A \subseteq \P(Z, \Psi)$, then
\begin{equation}
\label{EqInfScalar}
\forall \psi \in \Psi \colon \psi^\triup\of{\cl_\Psi \bigcup_{A \in \A} A} = \inf_{A \in \mathcal A} \psi^\triup(A)
\end{equation}
\end{lemma}

{\sc Proof.} Again, abbreviate $B = \cl_\Psi \bigcup_{A \in \A} A$. Since by (i) of Proposition \ref{PropPsiHull}
\[
\forall \psi \in \Psi, \; \forall A \in \A \colon \psi^\triup(B) \leq \psi^\triup\of{\bigcup_{A' \in \A} A'} \leq \psi^\triup(A),
\]
\begin{equation}
\label{ClosurScalar}
\forall \psi \in \Psi \colon \psi^\triup(B) \leq \psi^\triup\of{\bigcup_{A' \in \A} A'} \leq \inf_{A \in \A} \psi^\triup(A)
\end{equation}
follows; Lemma \ref{LemExtToLattice} (c) gives equality for the left inequality in \eqref{ClosurScalar}; and if ``$<$" would be true in the right inequality in \eqref{ClosurScalar}, then there would exist $\bar A \in \A$, $\bar a \in \bar A$ satisfying
\[
\psi(\bar a) < \inf_{A \in \A} \psi^\triup(A) \leq  \psi^\triup(\bar A),
\]
which is a contradiction. \pend

\medskip Equation \eqref{EqInfScalar} in Lemma \ref{LemInfScalar} can be written as
\begin{equation}
\label{EqInfStability}
\forall \psi \in \Psi \colon \inf_{A \in \A} \psi^\triup(A) = \psi^\triup\of{\inf_{A \in \A}A},
\end{equation}
and this crucial property is called {\em inf-stability} of the relation $\preceq_\Psi$ on $\P(Z, \Psi)$. The corresponding ``sup-stability" is in general not satisfied. The special case of a vector order $\preceq$ and support functions of closed convex sets as inf-extensions of continuous linear functionals can already be found in \cite[Lemma 4.14]{HamelEtAl15Incoll}. 

\begin{remark}
\label{RemSupExt}
Of course, one can start with the sup-extension of $\psi \in \Psi$ defined by $\psi^\trido(D) = \sup_{z \in D} \psi(z)$ and obtain a  corresponding hull operation as well as a set of closed sets which then has to be ordered by $\subseteq$. The theory is completely symmetric and adequate if maximization is the ultimate goal. 
\end{remark}

\begin{remark}
\label{RemEquivClasses}
The relation $\preceq_\Psi$ is not antisymmetric on $\P(Z)$ in general. Its symmetric part, i.e., the equivalence relation $\sim_\Psi$ defined by
\[
D \sim_\Psi E \quad \Leftrightarrow \quad D \preceq_\Psi E, \; E \preceq_\Psi D,
\]
can be expressed by $D \sim_\Psi E$ if, and only if, $\cl_\Psi D = \cl_\Psi E$. Thus, according to Proposition \ref{PropPsiLattice} the elements of $\P(Z, \Psi)$ can be understood as representatives of the equivalence classes with respect to $\sim_\Psi$. For particular $\Psi$, this has been observed, e.g., in \cite{HamelEtAl15Incoll}.
\end{remark}

\begin{remark}
\label{RemPointEmbedding}
The set $Z$ can be embedded into $\P(Z, \Psi)$ by defining a function $a \colon Z \to \P(Z, \Psi)$ by
\[
a(z) = \cl_\Psi\cb{z} = \bigcap_{\psi \in \Psi}\cb{y \in Z \mid \psi(z) \leq \psi(y)}
\]
whose values $a(z)$ include all $y \in Z$ with $z \preceq_\Psi y$, hence $a(z)$ is the upper level set of the restriction of $\preceq_\Psi$ at level $z$. By Proposition \ref{PropPsiHull} (i) $z \in a(z)$ and by (iii), $a(z) \in \P(Z, \Psi)$ for all $z \in Z$, and of course
\[
\cb{z_1} \preceq_\Psi \cb{z_2} \quad \Leftrightarrow \quad a(z_1) \supseteq a(z_2).
\]
Together with \eqref{EqExtOrder}, this gives an embedding of $(Z, \preceq)$ into $\of{\P(Z, \Psi), \supseteq}$. Equation \eqref{EqInfHull} implies $\psi(z) = \psi^\triup(a(z))$ for all $z \in Z$. 
\end{remark}

Very often, it is useful to reduce the number of elements in the set $\Psi$; for example, if $s > 0$ and $\Psi$ is replaced by $s\Psi = \cb{s\psi \mid \psi \in \Psi}$, then $\cl_\Psi D = \cl_{s\Psi} D$ for all $D \subseteq Z$ and  $\P(Z, \Psi) = \P(Z, s\Psi)$. This motivates the following definition.

\begin{definition}
\label{DefEquivFamilies}
Two sets $\Psi$ and $\Psi'$ of monotone functions on $(Z, \preceq)$ are called equivalent if
\[
\forall D \subseteq Z \colon \cl_\Psi D = \cl_{\Psi'} D.
\]
\end{definition}

If $\Psi$ and $\Psi'$  are equivalent, then \eqref{EqPsiPower} guarantees that $\P(Z, \Psi) = \P(Z, \Psi')$. One assumption which admits to introduce a ``smaller," but equivalent family reads as follows.

\begin{center}
\fbox{\; There is $\bar z \in Z$ satisfying $0 < \psi(\bar z) < +\infty$ for all $\psi \in \Psi$.\;} \qquad (A)
\end{center}

Indeed, if (A) is satisfied for the family $\Psi$, then $\Psi' = \cb{\frac{1}{\psi(\bar z)}\psi \mid \psi \in \Psi}$ is an equivalent family with $\psi'(\bar z) = 1$ for all $\psi' \in \Psi'$.

\medskip We shall give a few examples. The first one shows that every preorder can be represented by indicator functions of lower level sets. Here, the indicator functions in the sense of variational analysis is used, i.e., for $A \subseteq Z$ the function $I_A \colon Z \to \R\cup\cb{+\infty}$ is defined by $I_A(z) = 0$ if $z \in A$ and $I_A(z) = +\infty$ otherwise. The examples below also show that the union of $\Psi$-closed sets is not $\Psi$-closed in general since the topological closure and the convex hull turn out to be special cases.

\begin{example}
\label{ExIndicatorFam}
For the preordered set $(Z, \preceq)$, let $L(z) = \cb{y \in Z \mid y \preceq z}$ denote the lower and $U(z) = \cb{y \in Z \mid z \preceq y}$ the upper level set of the element $z \in Z$ with respect to $\preceq$, respectively. Then, the family of indicator functions $\mathcal I = \cb{I_{L(z)}}_{z \in Z}$ of the lower level sets with respect to $\preceq$ represents $\preceq$, i.e., for $z_1, z_2 \in Z$,
\[
z_1 \preceq z_2 \quad \Leftrightarrow \quad \forall z \in Z \colon I_{L(z)}(z_1) \leq I_{L(z)}(z_2).
\]
The proof is straightforward, so it is omitted. We will, however, verify that for $D \subseteq Z$,
\[
\cl_{\mathcal I} D = \bigcup_{z \in D}U(z) = \cb{y \in Z \mid \exists z \in D \colon z \preceq y} = \cb{y \in Z \mid D \cap L(y) \neq \emptyset}.
\]
Indeed, the definition of $\cl_{\mathcal I} D$ leads to
\[
\cl_{\mathcal I} D = \bigcap_{z \colon D \cap L(z) = \emptyset} \of{Z \bs L(z)}.
\]
It remains to show
\[
\bigcap_{z \colon D \cap L(z) = \emptyset} \of{Z \bs L(z)} = \cb{y \in Z \mid D \cap L(y) \neq \emptyset}.
\]

To show ``$\supseteq$": Take $y \in Z$ with $D \cap L(y) \neq \emptyset$. Assume $D \cap L(z) = \emptyset$ and $y \in L(z)$. Then there is $d \in D$ such that $d \preceq y \preceq z$, hence $d \in L(z)$ by transitivity, a contradiction. Hence $y \in Z \bs L(z)$.

To show ``$\subseteq$": Take $y \in \bigcap_{z \colon D \cap L(z) = \emptyset} \of{Z \bs L(z)}$. Assume $D \cap L(y) = \emptyset$. Then $y \in Z \bs L(y)$, i.e. $y \not\in L(y)$ which contradicts the reflexivity of $\preceq$.

Obviously, $\mathcal I$ does not satisfy assumption (A).
\end{example}

\begin{example}
As discussed, e.g., in \cite{EvrenOk11JME}, a multi-utility representation of the relation $\preceq$ (= preference in economic terms) on $Z$ is a family $\Psi = \mathcal U$ of functions $u \colon Z \to \R$ satisfying
\begin{equation}
\label{EqUtilityRep}
z_1 \preceq z_2 \quad \Leftrightarrow \quad \forall u \in \mathcal U \colon u(z_1) \leq u(z_2).
\end{equation}
It has been pointed out that each preorder has a (trivial) multi-utility representation (see \cite[Proposition 1]{EvrenOk11JME}) which is not very useful in most cases. See \cite{EvrenOk11JME}, the references therein and \cite{Evren14JET} for more relevant existence theorems, in particular in a topological setting. 

Since economists are usually interested in maximizing utility, the appropriate approach is via sup-extension: Defining
\[
u^\trido(D) = \sup_{z \in D} u(z) \quad \text{and} \quad \cl^\trido_{\mathcal U} D = \bigcap_{u \in \mathcal U}\cb{z \in Z \mid u(z) \leq u^\trido(D)}
\]
we obtain the complete lattice $\of{\mathcal P(Z, \mathcal U), \subseteq}$ where $\mathcal P(Z, \mathcal U) = \cb{D \subseteq Z \mid D = \cl^\trido_{\mathcal U} D}$, and the formulas for inf and sup are swapped: the infimum is an intersection, the supremum the $\mathcal U$-closure of the union.

Via \eqref{EqUtilityRep} and the preceding constructions, the utility maximization problem can be transformed into a set optimization problem for a complete lattice-valued function as follows. If $\mathcal Z_{ad} \subseteq Z$ is the set of admissible choices for the decision maker, solve
\[
\text{maximize} \quad b(z) \quad \text{subject to} \quad z \in \mathcal Z_{ad}
\]
with $b(z) = \cl^\trido_{\mathcal U} \cb{z}$. Because of \eqref{EqUtilityRep},
\[
z_1 \preceq z_2 \quad \Leftrightarrow \quad \cb{z_1} \preceq_\mathcal U \cb{z_2} \quad \Leftrightarrow \quad b(z_1) \supseteq b(z_2).
\] 
\end{example}

\begin{example}
Let $\Omega = \cb{\omega_1, \ldots, \omega_N}$ be a sample set and $\Pi$ a non-empty closed convex set of probability measures on $\Omega$. By $L^0$ we denote the linear space of all random variables $z \colon \Omega \to \R$, and it is clear that $L^0$ can be identified with $\R^N$. Let $u \colon \R \to \R\cup\cb{-\infty}$ be an increasing, concave function which we refer to as the utility function. By
\[
y \preceq z \quad \Leftrightarrow \quad \forall \pi \in \Pi \colon \E^\pi\sqb{u(y)} \leq \E^\pi\sqb{u(z)}
\]
a relation on $L^0$ is introduced which apparently is a preorder, but not a partial order in general. It is of course represented by the family
\[
\mathcal E = \cb{\E^\pi \circ u}_{\pi \in \Pi}
\]
where it is understood that $\of{\E^\pi \circ u}(z) := \E^\pi\sqb{u(z)} = -\infty$ if $P(\cb{\omega \in \Omega \mid u(y) = -\infty}) >0$. Relations of this type were introduced by Bewley \cite{Bewley86DP} (compare in particular \cite[Proposition 4]{GhirardatoEtAl03Econ} and \cite{RigottiShannon05Economet}) as one way to model incomplete preferences for which the incompleteness is due to uncertainty (= the ``true" probability measure is unknown). 

If there is $\bar r \in \R$ with $u(\bar r) > 0$, then $\mathcal E$ satisfies assumption (A). If $\dom u = \cb{r \in \R \mid -\infty < u(r)} \neq \emptyset$ (i.e., $u$ is a proper concave function), then this can always be arranged for by adding an appropriate constant to $u$, if necessary.

For $D \subseteq L^0$,
\[
\cl_{\mathcal E} D = \bigcap_{\pi \in \Pi}\cb{y \in L^0 \mid \E^\pi\sqb{u(y)} \leq \sup_{z \in D} \E^\pi\sqb{u(z)}},
\]
and
\[
b(x) = \cb{y \in L^0 \mid \forall \pi \in \Pi \colon  \E^\pi\sqb{u(y)} \leq \E^\pi\sqb{u(x)}} 
=  \cb{y \in L^0 \mid \sup_{\pi \in \Pi} \of{\E^\pi\sqb{u(y) - u(x)}} \leq 0} 
\]

The set-valued version of the utility maximization problem now is
\[
\text{maximize} \quad b(x) \quad \text{subject to} \quad x \in \mathcal X
\]
where $\mathcal X \subseteq L^0$ is the set of admissible alternatives; in many cases a ``budget set." 

The function $x \mapsto V(x) := \inf_{\pi \in \Pi} \of{\E^\pi\sqb{u(x) - u(y)}}$ was called, in a dynamic setting, the variational utility anchored at $y$ in \cite[Defintion~2.5]{DanaRiedel13JET}. The condition $V(x) \geq 0$ which appears in the definition of $b(x)$ above can be understood as sorting out all ``initial endowments" $y$ which are not better than $x$.
\end{example}

\begin{example}
\label{ExSSDviaAVaR}
Let $(\Omega, \F, P)$ be a probability space and $L^1 = L^1(\Omega, \F, P)$ the linear space of (equivalence classes of) integrable random variables over $(\Omega, \F, P)$. It is known (see \cite[Remark 4.49]{FoellmerSchied11}) that two random variables $x, y \in L^1$ are in relation with respect to second order stochastic dominance, denoted by $x \preceq_{SSD} y$, if, and only if,
\[
\forall \alpha \in (0, 1] \colon AV@R_\alpha(y) \leq AV@R_\alpha(x).
\]
The function $x \mapsto AV@R_\alpha(x)$ is called the Average Value at Risk (sometimes Conditional Value at Risk or Expected Shortfall, compare the discussion after \cite[Definition 4.48]{FoellmerSchied11}) and can be defined as follows:
\[
AV@R_\alpha(x) = \inf_{r \in \R}\cb{\frac{1}{\alpha}E[(r-x)^+] - r}.
\]
Defining the family of functions $\Psi = \cb{AV@R_\alpha}_{\alpha \in (0, 1]}$ on $Z = L^1$ we obtain another instance of an order relation defined via a family of real-valued functions. The question arises to find $\cl_\Psi$ and $\P(Z, \Psi)$.

Since $AV@R_\alpha(\One) = -1$ for all $\alpha \in (0, 1]$, assumption (A) can be satisfied. Here, $\One$ denotes the random variable which is equal to $1$ $P$-almost surely.
\end{example}

\begin{example}
\label{ExMultivariateFSD}
Again, let $(\Omega, \F, P)$ be a probability space and $L^0_d = L^0_d(\Omega, \F, P)$ the linear space of (equivalence classes of) random variables over $(\Omega, \F, P)$ with values in $\R^d$ for $d \geq 1$. Moreover, let $C \subseteq \R^d$ be a closed convex cone. The funcion $F_{x,C} \colon \R^d \to [0,1]$ defined by $F_{x,C}(z) = \inf\cb{P\cb{\omega \in \Omega \mid w^Tx(\omega) \leq w^Tz} \mid w \in C^+\bs\{0\}}$ is called the lower $C$-distribution function of $x$ with respect to $C$. If $d=1$ and $C = \R_+$, then $F_{x,\R_+}$ is just the (usual) cumulative distribution function of $x \in L^0:=L^0_1$, but for $d>1$ it is different from the joint distribution function even if $C = \R^d_+$ (see \cite{HamelKostner18JMVA}).

A random variable $y \in L^0_d$ is said to stochastically dominante $x \in L^0_d$ if
\[
\forall z \in \R^d \colon F_{x,C}(z) \leq F_{y,C}(z),
\]
and in this case we write $y \succeq^C_{FSD} x$. Defining a function $\psi_z \colon L^0_d \to [0,1]$ by $\psi_z(x) = F_{x,C}(z)$ one gets a representation of $\succeq^C_{FSD}$ via the family $\Psi = \cb{\psi_z}_{z \in \R^d}$.  If $d=1$ and $C = \R_+$, then this coincides with the usual first order stochastic dominance (see \cite[p. 93]{FoellmerSchied11}); for $d > 1$ this is a new concept proposed in \cite{HamelKostner18JMVA}.
\end{example}

In the remainder of this section, let $Z$ be a real linear space and $C \subseteq Z$ a convex cone with $0 \in C$. By 
\[
z_1 \leq_C z_2 \quad \Leftrightarrow \quad z_1 + C \supseteq z_2 + C
\]
a preorder on $Z$ is defined, i.e., $\leq_C$ is a reflexive and transitive relation on $Z$ which is compatible with the linear structure on $Z$.

\begin{example}
\label{ExSupportFunctions}
Let $Z$ be a separated locally convex, real linear space with topological dual $Z^*$ and $C \subseteq Z$ a convex cone with $0 \in C$ and $\cl C \neq Z$. Take $\Psi = C^+\bs\{0\}$ where
\[
C^+ = \cb{z^* \in Z^* \mid \forall z \in C \colon z^*(z) \geq 0}
\]
is the (positive) dual cone of $C$ including all continuous linear functionals on $Z$ which are monotone with respect to $\leq_C$. The above assumptions imply $C^+\neq\{0\}$. For $D \subseteq Z$, $\psi \in C^+$
\[
\psi^\triup(D) = \inf_{z \in D} \psi(z)
\]
is nothing, but the negative of the support function of $D$ taken at $-\psi$. By a separation argument one obtains for $D \in \P(Z)$
\[
\cl_{C^+\bs\{0\}} D = \cl\co(D + C),
\]
thus, with a slight abuse of notation,
\[
\P(Z, C^+) = \cb{D \in \P(Z) \mid D= \cl\co(D + C)}
\]
which is $\G(Z, C)$ in the notation of \cite{HamelEtAl15Incoll}. Proposition \ref{PropPsiLattice} yields that $(\G(Z, C), \supseteq)$ is a complete lattice with
\[
\inf \A = \cl\co \bigcup_{A \in \A} A \quad \text{and} \quad \sup \A = \bigcap_{A \in \A} A
\]
for $\A \subseteq \G(Z, C)$. The embedding function $a \colon Z \to \G(Z, C)$ is defined by $a(z) = z + \cl C$.
Finally, for $C = \cb{0}$ the class of all closed convex sets is obtained with $\Psi = C^+ = Z^*$. The complete lattice $(\G(Z, C), \supseteq)$ is the basic image space structure in set-valued convex analysis and optimization as established in \cite{Hamel09SVVAN, Hamel11Opt, Loehne11Book, HamelSchrage12JCA, HamelLoehne14JOTA} and recently surveyed in \cite{HamelEtAl15Incoll}.

If there is an element $\bar z \in C$ such that $z^*(\bar z) > 0$ for all $z^* \in C^+\bs\{0\}$, then the set $B^+ = \cb{z^* \in C^+ \mid z^*(\bar z) = 1}$ is a base of $C^+$ (see \cite[Theorem~2.2.12]{GoeRiaTamZal03Book}, applied to $C^+$ instead of $C$ with $(Z, Z^*)$ seen as a dual pair), and the two families $C^+\bs\{0\}$, $B^+$ are equivalent. This situation motivated Definition \ref{DefEquivFamilies} and Assumption (A).
\end{example}

\begin{example}
\label{ExTranslativeFunctions} Let $Z$ be a real linear space and $C$ a convex cone with $0 \in C$ generating the preorder $\leq_C$. Moreover, let $e \in C\bs(-C)$ be a fixed. Take $y \in Z$ and define a function $\tau_{y, e} \colon Z \to \R\cup\cb{\pm\infty}$ by
\[
\tau_{y, e}(z) = \inf\cb{t \in \R \mid y + te \in z + C} = \inf\cb{t \in \R \mid z - te \in y - C} 
\]
where $y - C = \cb{y - c \mid c \in C}$ is the usual Minkowski sum of $\cb{y}$ and $-C$. Then, the function $z \mapsto \tau_{y, e}(z)$ is monotone with respect to $\leq_C$ for each $y \in Z$ since for $z_1, z_2 \in Z$ with $z_1 \leq_C z_2$ we have $z_1 + C \subseteq z_2 + C$ and hence
\[
\tau_{y, e}(z_1) = \inf\cb{t \in \R \mid y + te \in z_1 + C} \leq \inf\cb{t \in \R \mid y + te \in z_2 + C} = \tau_{y, e}(z_2). 
\]
We set $\Psi = T(e) = \cb{\tau_{y, e}}_{y \in Z}$. If $y \in Z$ and $D \subseteq Z$, then
\[
\tau_{y, e}^\triup(D) = \inf_{z \in D} \tau_{y, e}(z) = \inf\cb{t \in \R \mid y + te \in D + C} =: \vp_{D+C, e}(y)
\]
where again $D + C$ is the usual Minkowski addition with the extension $\emptyset + C = \emptyset$.

A set $D \subseteq Z$ is called $e$-directionally closed  (see \cite{Schrage05Diploma, Hamel06R}) if
\[
\of{z \in Z \; \wedge \;  \cb{s_n}_{n \in \N} \subset \R \; \wedge \; \lim_{n \to \infty} s_n = 0 \;
\wedge \; \forall n \in \N \colon z + s_n e \in D} \; \Rightarrow \; z \in D.
\]
If $D \subseteq Z$ is not $e$-directionally closed, then one can add to $D$ all limits appearing in the last formula and obtains a set $\dc D$. This is a closure operation as also shown in \cite{Schrage05Diploma, Hamel06R}.  Moreover, it is known that
\[
\vp_{\dc(D+C), e}(y) = \vp_{D+C, e}(y) \leq 0 \quad \Leftrightarrow \quad y \in \dc(D+C).
\]
It can be shown that 
\[
\forall y \in Z \colon  \vp_{D+C, e}(y) \leq \tau_{y, e}(z) \quad \Leftrightarrow \quad z \in \dc(D+C).
\]
Indeed, if $z \in \dc(D+C)$, then $\vp_{D+C, e}(y) \leq \tau_{y, e}(z)$ follows from the definition of these two functions. Conversely, fix $z \in Z$ and assume
\[
\forall y \in Z \colon  \vp_{D+C, e}(y) \leq \tau_{y, e}(z).
\]
Then, in particular, $\vp_{D+C, e}(z) \leq \tau_{z, e}(z) = \inf\cb{t \in \R \mid z + te \in z + C} = \inf\cb{t \in \R \mid te \in C} = 0$ since $e \in C\bs(-C)$. Hence $z \in \dc(D+C)$ which proves the claim.

This shows
\[
\mathcal P(Z, \Psi) = \cb{A \in \P(Z) \mid A = \dc(A + C)}
\]
and $a(z) = \{z\} + \dc C$. Moreover, if $C$ is $e$-directionally closed, then $\Psi = \cb{\tau_{y, e}}_{y \in Z}$ is representing. If $Z$ is a topological linear space, $C \neq Z$ closed and $e \in \Int C$, then $\mathcal P(Z, \Psi) = \F(Z, C) := \cb{A \in \P(Z) \mid A = \cl(A + C)}$.
\end{example}

The following assumption can be used in order to reduce the set $\Psi = \cb{\tau_{y, e}}_{y \in Z}$ by transitioning to an equivalent family: there is $z' \in Z'$ where $Z'$ is the algebraic dual space of $Z$ such that $z'(e) = 1$ and $C \subseteq H^+(z') = \cb{z \in Z \mid z'(z) \geq 0}$. While the first part of this assumption can always be satisfied since $e \neq 0$ (in particular), sufficient conditions for the second part usually involve a separation argument with appropriate assumptions to $Z$ and $C$.

Under this assumption, let $\mathcal Y = \cb{y \in Z \mid z'(y) = -1}$ and $\bar z = e$. Then
\begin{align*}
\tau_{y, e}(\bar z) & = \inf\cb{t \in \R \mid y + te \in e + C}
	 = \inf\cb{t \in \R \mid y + te \in C} + 1 >  1
\end{align*}
for all $y \in \mathcal Y$ since $z'(y) + tz'(e) = -1 + t$ and $z'(C) \geq 0$ by assumption which means $y + te \in C$ can only be true for $t \geq 1$. Hence assumption (A) is satisfied.

Next, we claim that $\cb{\tau_{y, e}}_{y \in \mathcal Y}$ is representing. It is already clear from the above that $\tau_{y, e}$ is monotone w.r.t. $\leq_C$ even for all $y \in Z$. Next, assume
\[
\forall y \in \mathcal Y \colon \tau_{y, e}(z_1) \leq \tau_{y, e}(z_2)
\]
for $z_1, z_2 \in Z$. Then $\bar y := z_2 - (z'(z_2) + 1)e \in \mathcal Y$ and $\tau_{\bar y, e}(z_2) = z'(z_2) + 1$ (just check). Finally,
\begin{align*}
\tau_{\bar y, e}(z_1) & = \inf\cb{t \in \R \mid z_2 - z_1 + te \in C} + z'(z_2) + 1 \leq \tau_{\bar y, e}(z_2) = z'(z_2) + 1,
\end{align*}
hence $\varphi_{C,e}(z_2 - z_1) \leq 0$ which means $z_2-z_1 \in C$ since $C$ is assumed to be $e$-directionally closed.

\begin{example}
\label{ExOrientDist}
Let $Z$ be a normed space. For a set $A \subseteq Z$, its complement is $A^c = Z\bs A$. The oriented distance function $\Delta_A \colon Z \to \R\cup\cb{\pm\infty}$ for $A \subseteq Z$ is defined by
\[
\Delta_A(z) = d_A(z)-d_{A^c}(z)
\]
where $d_A(z) = \inf_{a \in A} \norm{a - z}$ and $d_\emptyset(z) = +\infty$. It was introduced in \cite[(1.5)]{HiriartUrruty79MOR} and is now widely used as a scalarization function in vector optimization, compare, for example, \cite{Zaffaroni03SIOPT}. The following statements are part of \cite[Prop 3.2]{Zaffaroni03SIOPT}. If $A \not\in \cb{Z, \emptyset}$, then

(i) $\Delta_A$ is real-valued and (globally) Lipschitz continuous with constant 1;

(ii) $\Delta_A(z)\leq 0$ if, and only if, $z \in \cl A$;

(iii) If $C \subseteq Z$ is a closed convex cone, then $\Delta_C$ is monotone with respect to $\leq_C$, i.e., $z_1 \leq_C z_2$ implies $\Delta_C(z_1) \leq \Delta_C(z_2)$.

The sup-extension of the function $\Delta_A$ was discussed in \cite{LiXu16Opt}, whereas its inf-extension can already be found in \cite{CrespiGinchevRocca06MMOR}. Here, we define the family 
\[
\Psi = \cb{p_y \colon Z \to \R \mid p_y(z) = \Delta_{z + C}(y), \, y \in Z}.
\]
Then $p_y(z) = \Delta_{z+C}(y) = \Delta_{C}(y-z)$ and $\Psi$ is representing for $\leq_C$ whenever $C$ is a closed convex cone. Indeed, while the monotonicity of the $p_y$'s follows from (iii) above, the converse implication follows by contradiction. Assume $p_{y}(z_1) \leq p_{y}(z_2)$ for all $y \in Z$ and $z_1 \not\leq_C z_2$, i.e., $z_2 - z_1 \not\in C$ for $z_2, z_1 \in Z$. Then, by (ii), $p_{z_2}(z_1) = \Delta_C(z_2 - z_1) > 0$ whereas $p_{z_2}(z_2) = \Delta_C(0) = 0$, hence $p_{z_2}(z_2) < p_{z_2}(z_1)$, a contradiction. It was already proven in \cite[Remark p.~84]{HiriartUrruty79MOR} that $\Psi$ is representing. 
 
The $\inf$-extension of $p_y \in \Psi$ is defined by
\[
\forall A \subseteq Z \colon p^\triup_y(A)=\inf_{z \in A}\Delta_{z+C}(y),
\]
and $p^\triup_y(A) = p^\triup_y(\cl(A+C))$ is immediate. Notice that $p^\triup_y(A) \neq \Delta_{A+C}(y)$. Indeed, let $Z=\R^2$, $C=\R^2_+$ and $A =\cb{(-1,0)^T,(0,-1)^T}+C$, $y = (1,1)^T$. Then $p^\triup_y(A) = -1$ while $\Delta_{A}(y) = -\sqrt{2}$.

The functions $p^\triup_y \colon \P(Z) \to \OLR$ have already been used, e.g., in \cite[Equation (3.15)]{LiuNgYang09MP}. For all $A \in \P(Z)$ it holds
\[
\cl_\Psi(A)=\cl(A + C).
\]
Indeed, as $p_y(A)=p_y(\cl(A+C))$ is true, we only need to prove the inclusion $\cl_\Psi(A) \subseteq \cl(A + C)$. 
Assume to the contrary that $z \in \cl_\Psi(A) \bs \cl(A + C)$ is true.  Then it exists $s>0$ such that $(z + \cb{z\in Z \mid \norm{z}\leq s})\cap \cl(A+C)=\emptyset$. But in this case, $p_z(A)\geq s > 0=p_z(\cl_\Psi A)$ which is a contradiction to Lemma \ref{LemExtToLattice}, \eqref{EqInfHull}. Hence $\P(Z, \Psi) = \F(Z, C)$ with the notation of \cite{HamelEtAl15Incoll}. One should note that this characterization of the complete lattice $(\F(Z, C), \supseteq)$ does not require any existence assumption as posed, for example, in \cite[Assumption 2.10]{AnsariChenYao17Opt} for the characterization of set order relations.
\end{example}

With Proposition \ref{PropPsiLattice} in view, we conclude this section by recalling an appropriate solution concept for complete lattice-valued minimization problems. The point of the following definition due to Heyde and L\"ohne \cite{HeydeLoehne11Opt} is that minimality and attainment of the infimum become two different concepts.

\begin{definition}
\label{DefLatticeSolution}
Let $X$ be a nonempty set, $(L, \leq)$ a complete lattice and $f \colon X \to L$ a function.

(a) A set $M \subseteq X$ is called a lattice-infimizer for $f$ if
\[
\inf_{x \in M} f(x) = \inf_{x \in X} f(x).
\]

(b) A point $\bar x \in X$ is called a lattice-minimizer for $f$ if
\[
x \in X, \; f(x) \leq f(\bar x) \quad \Rightarrow \quad f(x) = f(\bar x).
\]

(c) A set $M \subseteq X$ is called a lattice-solution of the problem
\[
\tag{P} \text{minimize} \quad f(x) \quad \text{over} \quad x \in X
\]
if $M$ is a lattice-infimizer and each $x \in M$ is a lattice-minimizer for $f$.

A lattice solution of (P) is called full if it includes all lattice-minimizers of $f$.
\end{definition}

One motivation for combining both of (a) and (b) into one solution concept is, of course, that $M = X$ trivially is an infimizer for $f$, another one that minimality in the sense of (b) is the most widely used minimality notion in vector and set optimization. This concept will mainly serve as a reference point for the new (approximate) optimality notions which are the subject of the next section. See the survey \cite{HamelEtAl15Incoll} for more details and references.

\section{Approximate solutions via families of scalar functions}
\label{SecSolConcepts}

In \cite[Definition 5.3]{HamelSchrage14PJO}, a solution concept for set optimization problems was introduced which is a modification of the one from \cite{HeydeLoehne11Opt, Loehne11Book} and uses scalarizations of $\G(Z, C)$-valued functions. In the following, this concept is extended in several directions. First, instead of solutions, approximate solutions will be defined. Secondly, the image space is not necessarily based on linear space constructions, but rather on the new approach presented in the previous section.

As before, $(Z, \preceq)$ is a pre-ordered set and  $\Psi$ a family of monotone (with respect to $\preceq$) functions $\psi \colon Z \to \OLR$. The set $\P(Z, \Psi)$ is defined through \eqref{EqPsiInfExt} and \eqref{EqPsiPower}. 

\begin{definition}
\label{DefApproxInf}
Let $\eps \geq 0$. A non-empty set $M \subseteq X$ is called an $(\eps, \Psi)$-infimizer of the function $f \colon X \to \P(Z, \Psi)$ if
\[
\forall \psi \in \Psi \colon 
\psi^\triup\of{\inf f[M]} \leq 
	\left\{
	\begin{array}{ccc}
	\psi^\triup\of{\inf f[X]} + \eps & : & \psi^\triup\of{\inf f[X]} \neq -\infty \\[.2cm]
	-\frac{1}{\eps} & : & \psi^\triup\of{\inf f[X]} = -\infty
	\end{array}
	\right.
\]
with the convention $-\frac{1}{0} = -\infty$. A $(0, \Psi)$-infimizer for $f$ is called $\Psi$-infimizer for $f$.
\end{definition}

In case of $\eps = 0$ the inequality in Definition \ref{DefApproxInf} can of course be replaced by $\psi^\triup\of{\inf f[M]} = \psi^\triup\of{\inf f[X]}$, and a lattice-infimizer always is a $(0, \Psi)$-infimizer for $f$.

As a stand-alone, the concept of $(\eps, \Psi)$-infimizers is not very interesting since the set $\dom f$ always is an $(\eps, \Psi)$-infimizer for all $\eps \geq 0$. Following the idea of \cite{HeydeLoehne11Opt} we complement the infimizer condition by an appropriate minimality property for its elements. Different from \cite{HeydeLoehne11Opt}, the family $\Psi$ is used to introduce alternatives to lattice-infimizers and -minimizers.

\begin{definition}
\label{DefApproxPsiMinPlus}
Let $\eps \geq 0$. An element $x_\eps \in X$ is called an $(\eps+, \Psi)$-minimizer of $f \colon X \to \P(Z, \Psi)$ if
\[
\forall \delta > 0, \;\exists \psi \in \Psi \colon 
\psi^\triup\of{f(x_\eps)} \leq 
	\left\{
	\begin{array}{ccc}
	\psi^\triup\of{\inf f[X]} + (\eps + \delta) & : & \psi^\triup\of{\inf f[X]} \neq -\infty \\
	-\frac{1}{\eps + \delta} & : & \psi^\triup\of{\inf f[X]} = -\infty,
	\end{array}
	\right.
\]
and the set of all $(\eps+, \Psi)$-minimizers of $f$  is denoted by $\Min(f,\eps+,\Psi)$.

An element $x_\eps \in X$ is called an $(\eps, \Psi)$-minimizer of $f \colon X \to \P(Z, \Psi)$ if
\[
\exists \psi \in \Psi \colon 
\psi^\triup\of{f(x_\eps)} \leq 
	\left\{
	\begin{array}{ccc}
	\psi^\triup\of{\inf f[X]} + \eps & : & \psi^\triup\of{\inf f[X]} \neq -\infty \\
	-\frac{1}{\eps} & : & \psi^\triup\of{\inf f[X]} = -\infty,
	\end{array}
	\right.
\]
and the set of all $(\eps, \Psi)$-minimizers of $f$ is denoted by $\Min(f,\eps,\Psi)$.
\end{definition}

Clearly, $\Min(f,\eps_1,\Psi) \subseteq \Min(f,\eps_2,\Psi)$ and $\Min(f,\eps_1+,\Psi) \subseteq \Min(f,\eps_2+,\Psi)$ for $0 \leq \eps_1 \leq \eps_2$.

\begin{remark}
\label{RemSolLinks}
(1) Obviously, $\Min(f,\eps,\Psi) \subseteq \Min(f,\eps+,\Psi)$ for all $\eps \geq 0$. Example \ref{ExEpsPlusVsEps} below shows that the inclusion can be strict.

(2) By definition,
\begin{equation}
\label{EqPlusDeltaMin}
\Min(f,\eps+,\Psi) = \bigcap\limits_{\delta>0}\Min(f,\eps + \delta,\Psi)
\end{equation}
for all $\eps \geq 0$. In particular, $\Min(f, 0+,\Psi) = \bigcap_{\eps>0}\Min(f, \eps,\Psi)$.

(3) From (1) and (2) one gets
\begin{align*}
\bigcap\limits_{\eps>0}\Min(f, \eps,\Psi) & \subseteq \bigcap\limits_{\eps>0}\Min(f, \eps+,\Psi) 
	= \bigcap\limits_{\eps>0} \bigcap\limits_{\delta>0}\Min(f,\eps+\delta,\Psi) =  \bigcap\limits_{\eps>0} \Min(f,\eps,\Psi),
\end{align*}
hence the inlusion is an equality. Therefore,
\begin{equation}
\label{EqIntersectEpsMin}
\Min(f,0+,\Psi) = \bigcap_{\eps>0}\Min(f, \eps,\Psi) = \bigcap\limits_{\eps>0}\Min(f, \eps+,\Psi).
\end{equation}
\end{remark}

\begin{remark}
\label{RemWeakMinimizers} Within the setting of Example \ref{ExSupportFunctions}, let $F \colon X \to Z\cup\{+\infty\}$ be a convex $C$-function which means that the set $\{(x,z) \in X \times Z \mid z \in F(x) + C\} $ is convex. Assume further that $\Int C \neq \emptyset$.  Let $f(x) = F(x) + C$ whenever $F(x) \neq +\infty$ and $f(x) = \emptyset$ otherwise. Then $\bar x \in X$ is a $(0, C^+)$-minimizer of $f$ if, and only if, it is a weakly minimal point of $F$. In a more general situation (even if $\Int C = \emptyset$), $(0, C^+)$-minimizers of set-valued functions $f$ were called $z^*$-solutions in \cite{HamelSchrage14PJO}.
\end{remark}

\begin{remark}
\label{RemVectorMin}
The first minimality notion for set-valued optimization problems (see e.g. \cite{Corley87JOTA, Luc89Book}) reads as follows. A point $(\bar x , \bar z) \in X \times Z$ is called a vector minimizer of the function $f \colon X \to \P(Z)$ if $\bar z \in f(\bar x)$ and $\bar z$ is a $\leq_C$-minimal point of $\bigcup\limits_{x \in X}f(x)$, i.e., $(\bar z - C)\cap(\bigcup\limits_{x\in X}f(x)) \subseteq \cb{\bar z} + C$. 

Let $X=\R$, $Z=\R^2$, $C=\R^2_+$, $\Psi=C^+ = \R^2_+$ and $f(x) = \cb{(s, \frac{1}{xs}) \mid s>0} + \R^2_+$ for all $x>0$ and $f(x)=\emptyset$ for $x \leq 0$. Then $\bigcup\limits_{x>0} f(x)=\Int \R^2_+$, and a vector minimizer does not exist. On the other hand, $\inf\cb{(1,0)^Tz \mid z \in f(x)}=0$ for all $x>0$, so each $x>0$ is a $(0,\Psi)$-minimizer.

Conversely, vector minimizers do not need to be $(0, \Psi)$-minimizers. For example, the function $f \colon \R \to \P(\R^2, \R^2_+)$ defined by
\[
f(x) = \left\{
	\begin{array}{ccc}
	\of{\begin{array}{cc}
	3x+3\\
	-x+3
		\end{array}} + \R^2_+
	 & : & x \in [-1,0] \\[.2cm]
	\of{\begin{array}{cc}
	x+3\\
	-3x+3
		\end{array}} + \R^2_+
	 & : & x \in [0,1] \\[.2cm]
	 \emptyset & : & \abs{x} > 1
	\end{array}
	\right.
\]
has plenty of vector minimizers which are not $(0, \R^2_+)$-minimizers.

While the issue of the first example is the missing infimum attainment for the scalarized problems, the crux of the second lies in the fact that the function $f$ (a nonconvex one) and the family $\Psi = C^+ = \R^2_+$ do not fit together: linear scalarizations mainly work well for convex problems.
\end{remark}

\begin{remark}
A relaxation of vector minimizers reads as follows (compare \cite{Kutateladze79SMD}). Let $\eps \geq 0$ and $h \in C\bs\{0\}$. A point $(x_\eps, z_\eps)$ is called an $\eps h$-minimizer of $F \colon X \to Z$ if
\[
\of{z_\eps - \eps h - C\bs\{0\}} \bigcap \of{\bigcup\limits_{x \in X} F(x)} = \emptyset.
\]

(i) The first example in Remark \ref{RemVectorMin} can easily be modified that there are $(\eps, \Psi)$-minimizers which are not $\eps h$-minimizers.

(ii) Conversely, the first component of an $\eps h$-minimizer does not need to be an $(\eps, \Psi)$-minimizer. Indeed, take $X = \R$, $Z = \R^2$, $C = H^+(e) = \cb{z \in \R^2 \mid h^Tz = z_1 + z_2 \geq {\f 0}}$, $\preceq = \leq_C$ with $e = (1, 1)^T$ and consider the function
\[
f(x) = \left\{
	\begin{array}{ccc}
	\cl\cb{z \in \R^2 \mid z_1 > 0, \; z_2 \geq \frac{x^2}{z_1}} & : & x \geq 0 \\[.2cm]
	\emptyset & : & x < 0
	\end{array}
	\right. 
\]
Clearly, $\cup_{x \in \R}f(x) = \Int \R^2_+$. For $x>0$, the inclusion $xe \in f(x)$ is true. If $0 < x \leq \eps$, then
\[
xe - \eps e - H^+(e)\bs\{0\} \bigcap \Int \R^2_+ = \emptyset,
\]
so every $x \in (0, \eps]$ is an $\eps e$-minimizer. On the other hand, with $\Psi = \cb{w \in C^+ \mid w_1 + w_2 = 2} = \cb{e}$ we get for $x > 0$
\[
\inf_{z \in f(x)} e^Tz = \inf_{s > 0}\cb{s + \frac{x^2}{s}} = 2x,
\]
and this means that all $x \in (\frac{\eps}{2}, \eps]$ are $\eps e$-minimizers, but not $(\eps, \Psi)$-minimizers since $\inf_{x>0} \inf_{z \in f(x)} = \inf_{x>0, \, s>0} \cb{s + \frac{x^2}{s}} = \inf_{x>0} 2x = 0$.

In \cite[Definition 3]{GreckschEtAl03Opt}, a relaxation of Kutateladze's definition was given which was upgraded to set-valued problems in \cite[Definition 17]{AlonsoDRodriguezM12JCAM}--the difference to $\eps h$-minimizers is that $\eps h$ is replaced by an element $h_\eps$ with $\norm{h_\eps} < \eps$ (in normed spaces). This concept is also different from $(\eps, \Psi)$-minimizers which can be shown by similar examples as above.
\end{remark}

Finally, we combine infimality and minimality into an approximate solution concept for set optimization problems.

\begin{definition}
\label{DefApproxPsiSolPlus}
A non-empty set $M \subseteq X$ is called an $(\eps+, \Psi)$-solution (an $(\eps, \Psi)$-solution) of the problem
\[
\tag{P} \text{minimize} \quad f(x) \quad \text{over} \quad x \in X
\]
if it is an $(\eps, \Psi)$-infimizer and each element of $M$ is an $(\eps+, \Psi)$-minimizer (an $(\eps, \Psi)$-minimizer) of $f$.
A $(0, \Psi)$-solution is called a $\Psi$-solution of (P). The set of all $(\eps+,\Psi)$-solutions and $(\eps,\Psi)$-solutions of (P) is denoted by $\Sol(f, \eps+,\Psi) \subseteq \P(X)$ and $\Sol(f, \eps,\Psi) \subseteq \P(X)$, respectively.
\end{definition}

A set $M \subseteq X$ is an $(\eps+,\Psi)$-solution (an $(\eps, \Psi)$-solution) of (P) if, and only if, $M \subseteq \Min(f, \eps+,\Psi)$ ($M \subseteq\Min(f,\eps,\Psi)$) and $\psi^\triup(f\sqb{M}) \leq \psi^\triup(f\sqb{X}) + \eps$ is true for all $\psi \in \Psi$, i.e., $M$ is an $(\eps,\Psi)$-infimizer consisting of only $(\eps+,\Psi)$-minimizers ($(\eps, \Psi)$-minimizers). 

Notably, $\Sol(f,\eps,\Psi)\subseteq \Sol(f,\eps+,\Psi)$, and $\Sol(f,\eps+,\Psi)$ is nonempty, if and only if, $\Min(f, \eps+, \Psi)$ is an element of $\Sol(f,\eps+,\Psi)$. It holds $\Sol(f,\eps+,\Psi)\subseteq \P(\Min(f,\eps+,\Psi))$ for all $\eps\geq 0$, and likewise for $\Sol(f,\eps,\Psi)$.

\begin{remark}
The novelty of Definition \ref{DefApproxPsiSolPlus} is twofold. First, approximate solutions are subsets of the preimage space rather than single points in the preimage space or of the graph of $f$. This feature--for lattice solutions--parallels Definition \ref{DefLatticeSolution}. Secondly, the values of approximate minimizers are not only compared to other function values, but also to the infimum of $f$ (see Definition \ref{DefApproxInf}). This key new feature sets Definition \ref{DefApproxPsiSolPlus} apart from previous concepts.

Scalarization procedures for set optimization problems, for example the one introduced in \cite{HamelLoehne06JNCA} which was subsequently used, for example, in \cite{HernandezRodriguezM07JMAA, GutierrezEtAl12NA, LiTeoZhang09NA}, also only take into account minimal elements with respect to set relations, but not the infimum. Therefore, similar examples as in Remark \ref{RemVectorMin} could be given which is not done here since it does not add more insights. 
\end{remark}

As in the scalar case, an $(\eps, \Psi)$-solution always exists {\f for $\eps >0$}. Note that--as in the scalar case--there is no compactness assumption necessary as used, e.g., in \cite[Corollary 25-27]{AlonsoDRodriguezM12JCAM}.

\begin{proposition}
\label{PropApproxSolExist}
Let $\eps > 0$ and $f \colon X \to \P(Z,\Psi)$ be a function. Then, there is an $(\eps, \Psi)$-solution of $\of{P}$.
\end{proposition}

{\sc Proof.} Fix $\psi \in \Psi$. By \eqref{EqInfStability}
\[
\psi^\triup\of{\inf f[X]} = \inf_{x \in X} \psi^\triup\of{f(x)}.
\]
If $\psi^\triup\of{\inf f[X]} \neq -\infty$, then there is $x_{\eps, \psi} \in X$ satisfying
\[
\psi^\triup\of{f(x_{\eps, \psi})} \leq \psi^\triup\of{\inf f[X]} + \eps,
\]
and if $\psi^\triup\of{\inf f[X]} = -\infty$, then there is $x_{\eps, \psi} \in X$ satisfying
\[
\psi^\triup\of{f(x_{\eps, \psi})} \leq -\frac{1}{\eps}.
\]
Define the set
\begin{multline*}
M = \sqb{\bigcup\limits_{\substack{\psi \in \Psi \\ \psi^\triup\of{\inf f[X]} \neq -\infty}}\negthickspace\negthickspace\negthickspace\cb{x \in X \mid \psi^\triup\of{f(x)} \leq \psi^\triup\of{\inf f[X]} + \eps}}
	\bigcup \\
	\sqb{\bigcup\limits_{\substack{\psi \in \Psi \\ \psi^\triup\of{\inf f[X]} = -\infty}}\negthickspace\negthickspace\negthickspace\cb{x \in X \mid \psi^\triup\of{f(x)} \leq -\frac{1}{\eps}}}.
\end{multline*}
This set $M$ is an $(\eps, \Psi)$-solution of $\of{P}$. Indeed, for a fixed  $\psi \in \Psi$ there is $x_{\eps, \psi} \in M$ such such (observe  \eqref{EqInfStability})
\[
\psi^\triup\of{\inf f[M]} = \inf_{x \in M} \psi^\triup\of{f(x)} \leq \psi^\triup\of{x_{\eps, \psi}} \leq
\left\{
\begin{array}{ccc}
 \psi^\triup\of{\inf f[X]} + \eps & : & \psi^\triup\of{\inf f[X]} \neq -\infty \\
 -\frac{1}{\eps} & : & \psi^\triup\of{\inf f[X]} = -\infty
\end{array}
\right.
\]
The second condition in Definition \ref{DefApproxPsiSolPlus} follows from the definition of $M$. \pend

\medskip As already discussed, some families $\Psi$ make more sense than others in particular situations. Thus, it becomes a task of the decision maker to carefully choose the set of potential scalarizations. The following two examples illustrate some difficulties which may occur. The first one shows that attention should be paid to the properness of the functions $\psi$ in Definition \ref{DefApproxPsiMinPlus}.

\begin{example}
\label{ExLexicographic}
Let $Z = \R^2$ and $C = \cb{z \in \R^2 \mid z_1 > 0, \;\text{or} \; z_1 = 0, z_2\geq 0}$, i.e., the lexicographic ordering cone in $\R^2$. With $e = (0,1)^T \in C\backslash (-C)$ we consider the family $\Psi= \cb{\tau_{y,e}}_{y \in \R^2}$ from Example \ref{ExTranslativeFunctions}. The set $\P(Z, \Psi)$ is the collection of sets $D \subseteq \R^2$ with $D = \dc D$. Especially, the function $\tau_{0, e}(z) = \inf\cb{t \in \R \mid -z + te \in -C}$ is an element of $\Psi$ with $\dom \tau_{0, e} = \cb{z \in \R^2 \mid z_1 \leq 0}$. Define a function $f \colon \R \to \P(Z, \Psi)$ by
\[
f(x) = (1+x, 1+x)^T + \cl C = (1+x, 1+x)^T + \cb{z \in \R^2 \mid z_1 \geq 0}
\]
for $x \geq 0$ and $f(x) = \emptyset$ whenever $x < 0$. Then,  in $(\P(Z, \Psi), \supseteq)$, $\inf f = (1,1)^T + \cb{z \in \R^2 \mid z_1 \geq 0}$ and $\inf f \cap \dom \tau_{0, e} = \emptyset$. This implies $\tau^\triup_{0, e}(\inf f) = +\infty$ which makes every $x$ an $(\eps, \Psi)$-minimizer. 
\end{example}

\begin{example}
\label{ExNonSolidCone}
Let $Z = \R^2$, $C = \cone{(1,1)^T} \cup\cb{(0,0)^T}$ and choose $\Psi = C^+\bs\{0\} = \cb{w \in \R^2\bs\{0\} \mid w_1 + w_2 \geq 0}$. Let $f \colon \R \to \P(\R^2, \Psi) = \G(\R^2, C)$ be given by $f(x) = (0, x)^T+C$ for $x \in \R$. Then, all values of $f$ are not comparable with each other (and hence each $x \in \R$ is a minimizer) and $\psi^\triup_w(f(x)) := \inf_{z \in f(x)}w^Tz = -\infty$ whenever $w$ is not parallel to $(1,0)^T$ in which case $\psi^\triup_w(f(x)) \equiv 0$ which again makes every $x \in \R^2$ an $(\eps, \Psi)$-minimizer. This is not too surprising. However, if $g(x) = x(1, 1)^T+C$ for $x \geq 0$ and $g(x) = \emptyset$ for $x < 0$ is considered, then for $w = (1, -1)^T$ or $w = (-1, 1)^T$ we have $\psi^\triup_w(g(x)) := \inf_{z \in g(x)} w^Tz \equiv 0$ is true, so again, every $x \geq 0$ is an $(\eps, \Psi)$-minimizer for each $\eps \geq 0$ which now is counterintuitive since $g(0)$ is the inifimum and $\bar x = 0$ is the only lattice-minimizer.
\end{example}

The following example shows that there can be more $(0+, \Psi)$-solutions than $(0, \Psi)$-solutions. The same function can be used to show that an $(\eps+, \Psi)$-minimizer is not an $(\eps, \Psi)$-minimizer, in general, for $\eps >0$.

\begin{example}
\label{ExEpsPlusVsEps}
Take $X = \R$, $Z = \R^2$, $C = \R^2_+$, $\preceq = \leq_C$ and consider the function
\[
f(x) = \left\{
	\begin{array}{ccc}
	\cl\cb{z \in \R^2 \mid z_1 > 0, \; z_2 \geq \frac{x^2}{z_1}} & : & x \geq 0 \\[.2cm]
	\emptyset & : & x < 0
	\end{array}
	\right. 
\]
Define $W = \cb{w \in \R^2_+ \mid w_1, w_2 > 0, \; w_1 + w_2 = 1}$ and take $\Psi = \cb{\psi_w}_{w \in W}$ with $\psi_w \colon \R^2 \to \R$ defined by $\psi_w(z) =  w^Tz$. Then $\psi^\triup_w(A) = \inf_{a \in A}w^Ta$ for $A \subseteq \R^2$. The function $f$ maps into $\P(\R^2, \Psi) \subsetneq \G(\R^2, \R^2_+)$ and $\inf f = \R^2_+ = f(0)$ in $(\P(\R^2, \Psi), \supseteq)$.

Since $\psi_w^\triup(\inf f) = 0$ for all $w \in W$, $M = \cb{\bar x}$ with $\bar x = 0$ is a $\Psi$-solution, and hence an $(\eps, \Psi)$-solution as well as an $(\eps+, \Psi)$-solution for all $\eps \geq 0$. On the other hand, since for $x > 0$ and $w \in W$
\[
\psi^\triup_w(f(x)) = \inf_{z_1 > 0}\cb{w_1z_1 + (1-w_1)\frac{x}{z_1}} = 2\sqrt{xw_1(1-w_1)} > 0,
\]
for each $\delta > 0$ there is $w \in W$ such that 
\[
0 < \psi^\triup_w(f(x)) \leq \delta.
\]
This shows that each $x \geq 0$ is a $(0+, \Psi)$-solution, but only $\bar x = 0$ also is a $(0, \Psi)$-solution.

One may observe that every $x \geq 0$ is a $(0, \Psi)$-solution for $\Psi = \cb{\psi_w}_{w \in \cl W}$ in which case $\P(\R^2, \Psi) = \G(\R^2, \R^2_+)$, and $f$ has the same infimum.
\end{example}

The following two examples show that lattice-minimality and $\Psi$-minimality are two different concepts.

\begin{example}[Frank's Example]\label{ExFrank} The following example is due to F. Heyde and was discussed in detail in \cite{HamelSchrage14PJO}. Take $X = Z = \R^2$, $C = \R^2_+$, $\preceq = \leq_C$ and consider the function $f \colon \R^2 \to \G(Z,C)$ defined by
\begin{align*}
f(x)=\left\{
         \begin{array}{ccc}
         \emptyset & : & x_1 < 0 \\
         \cb{z \in \R^2 \mid z_1\geq -x_1+x_2,\; z_2 \geq -x_1-x_2,\; z_1+z_2 \geq x_1} & : & x_1 \geq 0.
         \end{array}
       \right.
\end{align*}
A little sketching shows that every $x \geq 0$ is a lattice-minimizer in $(\G(Z,C), \supseteq)$, and $\inf_{x \in \R^2} f(x) = \cb{z \in \R^2 \mid z_1 + z_2 \geq 0}$.

Take $W = \cb{w \in \R^2_+ \mid w_1 + w_2 = 1}$ and $\Psi = \cb{\psi_w}_{w \in W}$ using the notation of Example \ref{ExEpsPlusVsEps}. Then $\P(Z, \Psi) = \G(Z,C)$ and
\[
\psi^\triup_w\of{\inf f} = 
	\left\{
      \begin{array}{ccc}
	-\infty & : & w \neq (1, 1)^T \\
	0 & : & w = (1, 1)^T.
	\end{array}
      \right.
\]
If $w \neq (1,1)^T$, then $\psi^\triup_w(f(x)) \in \R$ for all $x \in \dom f$, so 
\[
\forall \delta \geq 0 \colon \psi^\triup_w(f(x)) \leq -\frac{1}{\delta}
\]
can never be satisfied for an $x \in \dom f$. On the other hand, if $w = (1,1)^T$, then
\[
\forall x \geq 0 \colon \psi^\triup_w(f(x)) = x_1,
\]
hence $x \in \dom f$ is a $(0+, \Psi)$-minimizer (and a $(0, \Psi)$-minimizer at the same time) if, and only if, $x_1 = 0$.  This shows that lattice-minimizers can exist which are neither $(0+, \Psi)$- nor $(0, \Psi)$-minimizers.
\end{example}
 
\begin{example}
Take $X = \R$, $Z = \R^2$ and $C=\R^2_+$ with $\preceq = \leq_{\R^2_+}$. As in Example \ref{ExFrank}, take $\Psi = \cb{\psi_w}_{w \in W}$ with $W = \cb{w \in \R^2_+ \mid w_1 + w_2 = 1}$, a base of $C^+ = \R^2_+$. As before, $\P(Z, \Psi) = \G(Z, C)$. Let us consider the function $f \colon \R \to \G(\R^2, \R^2_+)$ defined by
\begin{align*}  
f(x)=\left\{
         \begin{array}{ccc}
          \emptyset & : & x < 0 \\
           \cb{z\in Z \mid  -x \leq z_1 \leq x,\; z_1 + z_2 \geq 0} + \R^2_+ & : & x \geq 0
         \end{array}
       \right.
\end{align*}

Then, $\inf_{x \in \R}f(x) = \cb{z \in Z \mid z_1+z_2 \geq 0}$, and each $x \geq 0$ is a $(0,\Psi)$-minimizer of $f$, but there is no minimizer with respect to $\supseteq$, i.e., no lattice-minimizer.

Moreover, $M \subseteq \R$ is a $(0,\Psi)$-solution of (P) if, and only if, $\sup M=+\infty$.
\end{example}

The following example shows that even if $\Psi \subseteq \bar\Psi$ is assumed, the sets of $(\eps,\Psi)$ and $(\eps,\bar\Psi)$ solutions are independent of each other.

\begin{example}
Let $X=\R$, $Z = \R^2$, $C = \R^2_+$ and $\lel=\leq_C$. Moreover, $\Psi=\cb{(1,0)^T, (0,1)^T}$ and $\overline\Psi=\Psi\cup\cb{(1,1)^T}$. Then, $a(z) = \bar a(z) = \cb{z} + \R^2_+$ for all $z \in\R^2$. Let $F \colon \R \to \R^2$ be defined by
\[
F(x)= \left\{
\begin{array}{ccc}
(x,1-x)^T & : &  0 \leq x \leq 1 \\
(x,0)^T & : &  1<x \\
(0,1 - x)^T & : & x < 0
\end{array}
\right.
\]
The $\P(\R^2, \Psi)$- and the $\P(\R^2, \bar\Psi)$-valued extensions of $F$ are equal and given by
\[
f(x)=\cb{F(x)}+ \R^2_+,
\]
with $\inf_{\Psi}f\sqb{\R} = \R^2_+$ and $\inf_{\bar\Psi}f\sqb{\R}= \R^2_+ \cap \cb{z \in \R^2 \mid z_1 + z_2 \geq 1}$.

The set of $(0,\Psi)$-minimizers of $f$ is $(-\infty, 0] \cup [1, \infty)$. A set $M \subseteq \R$ is a $(0, \Psi)$-infimizer if, and only if, $\cl M\cap\cb{x \in X \mid x\leq 0} \neq \emptyset$ and $\cl M\cap\cb{x\in X \mid x\geq 1} \neq \emptyset$. 

On the other hand, every $x \in \R$ is a $(0,\bar\Psi)$-minimizer of $f$, and a set $M \subseteq \R$ is a $(0,\bar\Psi)$ infimizer if, and only if, $\cl M\cap\cb{x \in \R \mid x \leq 0}\neq \emptyset$, $\cl M\cap\cb{x \in \R \mid x \geq 1}\neq \emptyset$ and $\cl M\cap\cb{x \in \R \mid 0 \leq x \leq 1}\neq \emptyset$.

Especially, $\cb{-1,2}$ is a $(0,\Psi)$-solution, but not a $(0,\bar\Psi)$-solution, while the open interval $\of{0,1}$ is a $(0,\bar\Psi)$-solution which is not a $(0,\Psi)$-solution.
\end{example}

The next two examples link approximate minimizers to approximate weakly efficient solutions in vector optimization. Note that the first one deals with the convex sitution (a convex, vector-valued function mapping into a locally convex space using linear scalarizations) while the second one also covers non-convex functions on topological linear spaces using translative scalarizations. For the same purpose, one could also use the oriented distance scalarization in the spirit of \cite{CrespiGinchevRocca06MMOR, Zaffaroni03SIOPT}.

If $X, Z$ are two linear spaces, $S \subseteq X$, $C \subseteq Z$ a convex cone and $F \colon X \to Z$ a function, then the function $f \colon X \to \P(Z)$ defined by $f(x) = \cb{F(x)} + C$ for $x \in S$ and $f(x) = \emptyset$ for $x \not\in S$ is called the (set-valued) $C$-extension of $F$.

\begin{example}
\label{ExWeakEffLinear}
Let $X$ be a linear space and $Z$ a locally convex Hausdoff topological space with dual $Z^*$, let $C \subseteq Z$ be a closed convex cone with $\Int C \neq \emptyset$ and fix $e \in \Int C$. Then, the set $B^+ = \cb{z^* \in C^+ \mid z^*(e) = 1}$ is a base of $C^+$. We choose $\Psi = B^+$ and denote by $\psi_{z^*} \in \Psi$ the element in $\Psi$ which coincides with $z^*$. In this case, $\P(Z, \Psi) = \G(Z, C)$ (see Example \ref{ExSupportFunctions}).

Take a set $\emptyset \neq S \subseteq X$ and a function $F \colon S \to Z$. The $C$-extension $f$ of $F$ actually maps into $\P(Z, \Psi)$ since $C$ is a closed convex cone. The set of approximate weakly efficient solutions for $\eps \geq 0$ of the vector minimization problem for $F$ is 
\[
\WEff_{\eps e}(F,S) = \cb{\bar x \in S \mid \forall x \in S \colon F(x) \not\in F(\bar x) - \eps e - \Int C}.
\]
If $S$ is convex and $F$ is $C$-convex, then
\[
\forall \eps \geq 0 \colon \WEff_{\eps e}(F,S) = \Min(f_C, \eps, B^+) =  \Min(f_C, \eps+, B^+).
\]
Indeed, if $\bar x \in \WEff_{\eps e}(F,S)$ then $F[S] \cap \of{F(\bar x) - \eps e - \Int C} = \emptyset$, so one can separate the two convex sets $F[S]$ and $F(\bar x) - \eps e - C$ getting $z^* \in B^+$ with $z^*(F(\bar x) - \eps e) \leq z^*(F(x))$ for all $x \in S$ which produces $z^*(F(\bar x)) \leq \inf_{x \in S}(z^* \circ F)(x) + \eps$. Hence $\bar x \in \Min(f, \eps, B^+) \subseteq  \Min(f, \eps+, B^+)$. Conversely, if $\bar x \in \Min(f, \eps+, B^+)\bs \WEff_{\eps e}(F,S)$, then there is $x \in S$ such that
\[
F(x) \in F(\bar x) - \eps e - \Int C
\] 
which yields $z^*(F(x) < z^*(F(\bar x)) - \eps$ for all $z^* \in B^+$ and therefore
\[
\forall z^* \in B^+ \colon \psi^\triup_{z^*}(f(\bar x)) > \psi^\triup_{z^*}(f(x)) + \eps \geq \inf_{x \in S}\psi^\triup_{z^*}(f(x)) + \eps = \psi^\triup_{z^*}(\inf_{x \in S}f(x)) + \eps
\]
where the last equation is inf-stability (see \eqref{EqInfStability}). This means that $\bar x$ cannot be an $(\eps+, B^+)$-minimizer for $f$ which contradicts the assumption. 
\end{example}

\begin{example}
\label{ExWeakEffTranslative}
Let $X$ be a linear space and $Z$ be a topological linear space, let $C \subsetneq Z$ be a closed convex cone with $\Int C \neq \emptyset$ and fix $e \in \Int C$. We choose $\Psi = \cb{\tau_{y, e}}_{y \in Y}$ with 
\[
\tau_{y, e}(z) = \inf\cb{t \in \R \mid y + te \in z + C} = \inf\cb{t \in \R \mid z - te \in y - C}.
\]
Under the above assumption, it follows, e.g., from Theorem 2.3.1 and Proposition 2.3.4 in \cite{GoeRiaTamZal03Book} that for all $y \in Z$,
$\tau_{y, e}$ is finite-valued, continuous, convex and strictly monotone, i.e.,
\[
z_2 - z_1 \in \Int C \quad \Rightarrow \quad \tau_{y, e}(z_1) < \tau_{y, e}(z_2).
\]
Then $\P(Z, \Psi) = \F(Z, C) = \cb{D \subseteq Z \mid D = \cl(D+C)}$ (see Example \ref{ExTranslativeFunctions}). Moreover, $\tau_{y, e}(z) < 0$ if, and only if, $z \in y - \Int C$.

Let $S \subseteq X$ and $F \colon S \to Z$ be a non-empty set and a (vector-valued) function, respectively. Then, the $C$-extension $f$ of $F$ maps into $\G(Z, C):= \cb{D \subseteq Z \mid D = \cl\co(D+C)} \subseteq \F(Z,C)$, and for $\eps \geq 0$ one has
\[
\WEff_{\eps e}(F,S) = \Min(f, \eps, \Psi).
\]
Indeed, assume first $\bar x \in \WEff_{\eps e}(F,S)$. Denote $\bar z = F(\bar x)$. Then $\sqb{F[S] + \eps e} \cap \sqb{F(\bar x) - \Int C} = \emptyset$ and hence $\tau_{\bar z, e}(F(x) + \eps e) = \tau_{\bar z, e}(F(x)) + \eps \geq 0$ for all $x \in S$. On the other hand, $\tau_{\bar z, e}(\bar z) = 0$. Altogether, $\tau_{\bar z, e}(\bar z) = 0 \leq \inf_{x \in S}\tau_{\bar z, e}(F(x)) + \eps$ which shows $\bar x \in \Min(f, \eps, \Psi)$. Conversely, assume $\bar x \in  \Min(f, \eps, \Psi)\bs \WEff_{\eps e}(F,S)$. Then, there are $x \in S$ and $y \in Z$ such that $F(x) +  \eps e \in  F(\bar x)  - \Int C$ and
\[
\forall x \in S \colon \tau_{y,e}(F(\bar x)) \leq \tau_{y,e}(F(x)) + \eps.
\]
Since $\tau_{y,e}$ is strictly monotone, one gets from these two relationships
\[
\tau_{y,e}(F(x) + \eps e) =  \tau_{y,e}(F(x)) + \eps <  \tau_{y,e}(F(\bar x)) \leq \tau_{y,e}(F(x)) + \eps,
\]
a contradiction.
\end{example}

To conclude this section, approximate solutions are discussed when the representing family is the one consisting of indicator functions.

\begin{example}
\label{ExIndicatorFamCont}
Let $f \colon X \to \mathcal P(Z, \mathcal I)$ be a function and consider the problem
\[
\text{minimize} \quad f(x) \quad \text{subject to} \quad x \in X.
\]
What are approximate solutions with respect to the representing family $\mathcal I = \cb{I_{L(z)}}_{z \in Z}$?

Let $\eps \geq 0$. Then $M_\eps$ is an $(\eps, \mathcal I)$-infimizer if
\[
\forall y \in Z \colon I^\triup_{L(y)}\of{\inf f\sqb{M_\eps}} \leq  I^\triup_{L(y)}\of{\inf f\sqb{X}} + \eps. 
\]
First, observe 
\[
I^\triup_{L(y)}(D) = \inf_{z \in D} I_{L(y)}(z) = 
	\left\{
	\begin{array}{ccc}
	0 & : & D \cap L(y) \neq \emptyset \\
	+\infty & : & D \cap L(y) = \emptyset
	\end{array}
	\right\}
	= I_D(y)
\]
for $D \in \mathcal P(Z, \mathcal I)$ and secondly, $\inf f\sqb{M_\eps} = \bigcup_{x \in M_\eps} f(x)$ since $f$ maps into $\mathcal P(Z, \mathcal I)$ (see Example \ref{ExIndicatorFam}). If $\bar z \in f\sqb{X}\bs f\sqb{M_\eps}$, then $\bar z \not\in f(x)$ for all $x \in M_\eps$, hence
\[
I^\triup_{L(y)}\of{\inf f\sqb{X}} = 0 \quad \text{and} \quad I^\triup_{L(y)}\of{\inf f\sqb{M_\eps}} = +\infty.
\]
So, $M_\eps$ is an $(\eps, \mathcal I)$-infimizer if, and only if, $\inf f\sqb{M_\eps} = \inf f\sqb{X}$, i.e., $M_\eps$ is an infimizer.

On the other hand, let $x_\eps \in X$ be an $(\eps, \mathcal I)$-minimizer of $f$, i.e.,
\[
\exists y \in Z \colon I^\triup_{L(y)}\of{f(x_\eps)} \leq I^\triup_{L(y)}\of{\inf f\sqb{X}} + \eps.
\]
Since $I^\triup_{L(y)}\of{f(x_\eps)}  = I_{f(x_\eps)}(y)$ and $I^\triup_{L(y)}\of{\inf f\sqb{X}} = I_{\inf f[X]}(y)$, the above inequality means
\[
I_{f(x_\eps)}(y) \leq I_{\inf f[X]}(y) + \eps
\]
which is obviously satisfied for $y = z_\eps \in f(x_\eps)$. So, $x_\eps$ is an $(\eps, \mathcal I)$-minimizer if, and only if, $f(x_\eps) \neq \emptyset$. Altogether, this shows that the family $\mathcal I$ is a very good choice if one is only interested in infimizers, but it is a bad one if one looks for (approximate) minimizers. However, this once again emphasizes the role of the infimum in set optimization.
\end{example}

\section{Existence theorems and well-posedness for set optimization}
\label{SecExWellPo}

Under compactness and lower semicontinuity assumptions, it will be shown that there exists a $(0, \Psi)$-solution of problem (P). The result can be seen as a Weierstrass-type existence theorem for set optimization problems. In contrast to \cite{HeydeLoehne11Opt}, the (lower) domination property does not play a role here, but rather the scalar (Weierstrass) extreme value theorem and Hausdorff convergence of approximate solutions.

In this section, it is assumed throughout that $X$ is a separated topological linear space over the reals. For the convenience of the reader, we state a Weierstrass type result for improper extended real-valued functions. We think  that there should be a reference for it, but could not find one.

\begin{proposition}
\label{PropImpropWeier}
Let  $\vp \colon X \to \R\cup\cb{\pm\infty}$ an extended real-valued function. If there is $\bar r \in \R$ such that the sublevel sets 
\[
L_\vp(r) = \cb{x \in X \mid \vp(x) \leq r}
\]
are compact for $\inf_{x \in X} \vp(x) < r \leq \bar r$, then there is $\bar x \in X$ with $\vp(\bar x) = \inf_{x \in X} \vp(x)$.
\end{proposition}

{\sc Proof.} Denote $\bar \vp = \inf_{x \in X} \vp(x)$. If $\bar \vp < r$ for $r \in \R$, then $L_\vp(r) \neq \emptyset$. Moreover, $L_\vp(r) \subseteq L_\vp(s)$ for $r \leq s$. Since for $\bar \vp < r \leq \bar r$ the sets $L_\vp(r)$ are compact, one has
\[
\bigcap_{r \in (\bar \vp, \bar r]} L_\vp(r) \neq \emptyset
\]
by Cantor's intersection theorem. Clearly, any $\bar x \in \cap_{r \in (\bar \vp, \bar r]} L_\vp(r)$ 
must satisfy $\vp(\bar x) = \inf_{x \in X} \vp(x)$. \pend

\medskip The above result applies in particular if $\inf_{x \in X} \vp(x) = -\infty$. The basic result reads as follows.

\begin{theorem}
\label{ThmPsiWeierstrass} 
Let $f \colon X \to \P(Z,\Psi)$ be such that for each $\psi \in \Psi$ the function $\psi^\triup \circ f \colon X \to \R\cup\cb{\pm\infty}$ is lower semicontinuous and there is $r_\psi \in \R$ such that $\inf(\psi^\triup \circ f)[X] < r_\psi$ and $L_{\psi^\triup \circ f}(r_\psi)$ is compact. Then, for each $\psi \in \Psi$ there is $x_\psi \in X$ such that 
\begin{equation}
\label{EqPsiCircFMinimizer}
(\psi^\triup \circ f)(x_\psi) = \inf(\psi^\triup \circ f)\sqb{X}.
\end{equation}
Moreover, the set $\cb{x_\psi \in X \mid x_\psi \; \text{satisfies \eqref{EqPsiCircFMinimizer}}, \; \psi \in \Psi}$ is a $(0, \Psi)$-solution of (P).
\end{theorem}

{\sc Proof.} Fix $\psi \in \Psi$. Since the sets $L_{\psi^\triup \circ f}(r)$ for $-\infty < r \leq r_\psi$ are closed and subsets of $L_{\psi^\triup \circ f}(r_\psi)$ they are also compact, and one can apply Proposition \ref{PropImpropWeier} to get $x_\psi \in X$ satisfying \eqref{EqPsiCircFMinimizer}. Moreover, $x_\psi \in \Min(f, 0, \Psi)$ by definition.  Since this argument is valid for each $\psi \in \Psi$, 
$\Min(f, 0, \Psi)$ is a $(0, \Psi)$-infimizer of $f$ and hence a $(0, \Psi)$-solution for (P). \pend

\begin{theorem}
\label{ThmBasicWeierstrass}
Let $f \colon X \to \P(Z,\Psi)$ and $\eps_0 > 0$ be such that $\Min(f,\eps+,\Psi)$ is compact for all $0 < \eps \leq \eps_0$. Then $\Min(f,0+,\Psi)$ is nonempty and compact, and $\Min(f,\eps+, \Psi)$ Hausdorff converges to $\Min(f, 0+,\Psi)$ as $\eps \to 0$.

If, additionally, $\psi^\triup\circ f$ is lower semicontinuous for all $\psi \in \Psi$, then 
$\Min(f,0,\Psi)$ is non-empty and a $(0,\Psi)$-solution of problem (P).
\end{theorem}

{\sc Proof.} According to $\eqref{EqIntersectEpsMin}$,
\[
\Min(f,0+,\Psi) = \bigcap\limits_{\eps>0}\Min(f,\eps+,\Psi),
\]
and this set is nonempty and compact by Cantor's intersection theorem.

Assume that $\Min(f,\eps+,\Psi)$ does not upper Hausdorff converge to $\Min(f,0+,\Psi)$ as $\eps$ converges to $0$. Then, there exists a neighborhood $U$ of $0 \in X$ such that 
\[
\forall \delta > 0, \; \exists \eps \in (0, \delta) \colon \Min(f,\eps+,\Psi)\bs\of{\Min(f, 0+, \Psi) + U} \neq \emptyset.
\]
Let $\cb{\eps_i}_{i \in I} \subseteq (0, \eps_0]$ be a decreasing net  with $\lim_{i \in I} \eps_i = 0$. For each $i \in I$ pick $x_i \in \Min(f,\eps_i+,\Psi)\bs(\Min(f, 0+, \Psi) + U)$. Then $\cb{x_i}_{i \in I} \subseteq \Min(f,\eps_0+,\Psi)$. Since the latter set is compact, there is a convergent subnet $\cb{x_{i_j}}_{j \in I}$ of $\cb{x_i}_{i \in I}$ which converges to $\bar x \in \Min(f,\eps_0+,\Psi)$. Fix $i \in I$. Then
\[
\forall i' \in I, i' \succ i \colon x_{i'} \in \Min(f,\eps_i+,\Psi),
\]
and the latter set is compact by assumption. Hence $\bar x \in \Min(f,\eps_i+,\Psi)$ for all $i \in I$ which yields
\[
\bar x \in \bigcap_{\eps > 0}\Min(f, \eps+, \Psi)=\Min(f, 0+, \Psi)
\]
since $\lim_{i \in I} \eps_i = 0$. Because $\cb{x_{i_j}}_{j \in I}$ converges to $\bar x$ there is $j_0 \in J$ such that
\[
\forall j \succ j_0 \colon x_{i_j} \in \bar x + U \subseteq \Min(f, 0+, \Psi) + U
\]
which contradicts the assumption. Therefore, $\Min(f,\eps+,\Psi)$ upper Hausdorff converge to $\Min(f,0+,\Psi)$.

Obviously, $\Min(f,\eps+,\Psi)$ lower Hausdorff converges to $\Min(f,0+,\Psi)$ since $\Min(f,0+,\Psi)\subseteq \Min(f,\eps+,\Psi)$ for all $\eps \geq 0$.

Finally, assume $\psi^\triup\circ f$ is lower semicontinuous for all $\psi\in \Psi$. Then, the level sets of $\psi^\triup \circ f$ are closed. 

If $\inf(\psi^\triup \circ f)[X] \neq -\infty$, then in particular the sets
\[
L_{\psi^\triup \circ f}(\inf\of{\psi^\triup \circ f}\sqb{X} + \eps) = \cb{x \in X \mid (\psi^\triup \circ f)(x) \leq \inf\of{\psi^\triup \circ f}\sqb{X} + \eps}
\]
for $\eps > 0$ are closed. Since $L_{\psi^\triup \circ f}(\inf\of{\psi^\triup \circ f}\sqb{X} + \eps) \subseteq \Min(f, \eps, \Psi) \subseteq \Min(f, \eps+, \Psi)$ and the latter set is compact, so is $L_{\psi^\triup \circ f}(\inf\of{\psi^\triup \circ f}\sqb{X} + \eps)$. If $\inf(\psi^\triup \circ f)[X] = -\infty$, the same argument applies to the sets $L_{\psi^\triup \circ f}(-\frac{1}{\eps}) = \cb{x \in X \mid (\psi^\triup \circ f)(x) \leq -\frac{1}{\eps}} \subseteq \Min(f, \eps, \Psi)$. Hence $\Min(f,0,\Psi)$ is non-empty and a $(0,\Psi)$-solution of problem (P) by Theorem \ref{ThmPsiWeierstrass}. \pend

\medskip Of course, the statement of Theorem \ref{ThmBasicWeierstrass} remains true if $\Min(f,\eps+,\Psi)$ in its assumptions is replaced by $\Min(f,\eps,\Psi)$ for $\eps>0$. Existence result for vector optimization problems such as Corollary \ref{CorVectWeakEff} below can be obtained (and even improved) as special cases of the above set optimization result. In order to give the set-up, the following definition is adapted from \cite[Definition 5.1]{Luc89Book}.

\begin{definition}
\label{DefCCont}
Let $Z$ be topological linear space, $C \subseteq Z$ a convex cone with $0 \in C$ and $S \subseteq X$ a non-empty set. A function $F \colon S \to Z$ is called $C$-continuous at $\bar x \in S$ if for each neighborhood $W$ of $0 \in Z$ there exists a neighborhood $U$ of $0 \in X$ such that
\[
\forall x \in (\bar x + U)\cap S \colon F(x) \in F(\bar x) + W + C.
\]
The function $F$ is called $C$-continuous on $S$ if it is $C$-continuous at every $x \in S$.
\end{definition}

\begin{corollary}
\label{CorVectWeakEff}
Within the setting of Definition \ref{DefCCont}, let $\Int C \neq \emptyset$ and $F \colon S \to Z$ be a $C$-continuous function. Moreover, assume that $\WEff_{\eps_0 e}(F,S)$ is compact relative to $S$ for some $\eps_0>0$. If either

(i) $e \in \Int C$ and $\Psi = T(e) = \cb{\tau_{y,e} \mid y \in Z}$ or

(ii) $Z$ is a separated locally convex space, $\Psi = B^+ = \cb{z^* \in C^+ \mid z^*(e) = 1}$ and $F$ is $C$-convex,

then $\WEff(F,S)$ is nonempty, compact relative to $S$, $\WEff_{\eps e}(F,S)$ upper Hausdorff converges to $\WEff(F,S)$ for $\eps \to 0$, and $\WEff(F,S)$ is a $(0,\Psi)$-solution for (P) where $f$ is the $C$-extension of $F$.
\end{corollary}

We precede the proof of the corollary with two auxiliary lemmas.

\begin{lemma}
\label{LemCContClosed}
If $F \colon S \to Z$ is a $C$-continuous function, then $\WEff_{\eps e}(F,S)$ is closed relative to $S$ for all $\eps \geq 0$.
\end{lemma}

{\sc Proof.} Take $\bar x \in \of{\cl \WEff_{\eps e}(F,S)} \cap S$. Then, for each neighborhood $U$ of $0 \in X$ one has
\begin{equation}
\label{EqCloseEff}
(\bar x + U) \cap \WEff_{\eps e}(F,S) \neq \emptyset.
\end{equation}
Assume that $\bar x \not\in \WEff_{\eps e}(F,S)$. Then, there are $\hat x \in S$ and $\delta > 0$ such that $F(\hat x) \in F(\bar x) - (\eps + \delta) e - \Int C \subseteq F(\bar x) - \eps e - \Int C$. Since $-\delta e + \Int C$ is a neighborhood of $0 \in Z$ and $F$ is $C$-continuous there is a neighborhood $U_\delta$ of $0 \in X$ such that 
\[
\forall x \in \bar x + U_\delta \colon F(x) \in F(\bar x) - \delta e + \Int C + C = F(\bar x) - \delta e + \Int C.
\]
Take $x_\delta \in \bar x + U_\delta$. Then,
\[
F(\hat x) \in F(\bar x) - (\eps + \delta)e - \Int C \subseteq F(x_\delta) + \delta e - \Int C - (\eps + \delta)e - \Int C = F(x_\delta) - \eps e - \Int C.
\]
This is a contradiction since one can choose $ x_\delta \in  \WEff_{\eps e}(F,S)$ by \eqref{EqCloseEff}. \pend

\begin{lemma}
\label{LemCContClosed}
If the assumptions of Corollary \ref{CorVectWeakEff} including (i) are satisfied, then the functions $\tau_{y,e} \circ F \colon X \to \R\cup\cb{\pm\infty}$ for $y \in Z$ are lower semicontinuous relative to $S$. The same holds in case of (ii) for the functions $z^* \circ F \colon X \to \R\cup\cb{\pm\infty}$ for $z^* \in B^+$.
\end{lemma}

{\sc Proof.} We give the proof for (i) and omit the very similar one for (ii). Fix $\bar x \in S$. Since $e \in \Int C$, the set $-\eps e + \Int C$ is a neighborhood of $0 \in Z$ for all $\eps > 0$. By $C$-continuity of $F$, there exists a neigborhood $U$ of $0 \in X$ such that
\[
\forall x \in (\bar x + U) \cap S \colon F(x) \in F(\bar x) -\eps e + \Int C + C \subseteq F(\bar x) -\eps e + C.
\]
Since the functions $\tau_{y, e}$ are monotone (see Example \ref{ExTranslativeFunctions})
\[
\forall x \in (\bar x + U) \cap S \colon \tau_{y, e}(F(\bar x) -\eps e) = \tau_{y, e}(F(\bar x))  - \eps \leq F(x)
\] 
follows which means that $x \mapsto (\tau_{y, e} \circ F)(x) = \tau_{y, e}(F(x))$ is lower semicontinuos at $\bar x \in S$.
\pend

\medskip {\sc Proof of Corollary \ref{CorVectWeakEff}.} Case (i): From Example \ref{ExWeakEffTranslative}, one gets $\WEff_{\eps e}(F,S) = \Min(f, \eps, T(e))$ for all $\eps \geq 0$. From Theorem \ref{ThmBasicWeierstrass} one can conclude that $\WEff(F,S) = \Min(f, 0, T(e))$ is non-empty, compact and that $\WEff_{\eps e}(F,S) = \Min(f, \eps, T(e))$ Hausdorff converges to it. The rest follows from the second part of Theorem \ref{ThmBasicWeierstrass} and Lemma \ref{LemCContClosed}. Case (ii) follows similarly. \pend

\medskip Note that Corollary \ref{CorVectWeakEff} above differs from \cite[Theorem 4.2]{CrespiGuerraggioRocca07JOTA}  since here it is assumed that $\WEff_{\eps_0 e}(F,S)$ to be compact, rather then the sets $\WEff_{\eps e}(F,\lev_y)$ for all nonempty level sets $\lev_y = \cb{x \in X \mid F(x) \in y-C}$, and $C$-quasiconvexity is not assumed.

In the setting of Corollary \ref{CorVectWeakEff}, $F$ is well-posed in the sense of \cite[Definition 3.4]{CrespiGuerraggioRocca07JOTA}. We conclude this section by introducing a notion of well-posedness for set optimization problems which fits into the framework of this paper.

\begin{definition}
Let $f \colon X \to \P(Z,\Psi)$ be given. Then $\cb{M_i}_{i\in I}$ is called a $\Psi$-minimizing net, if for all $\eps>0$ there is a $i_\eps\in I$ such that for all $i>i_\eps$ the set $M_i$ is a $(\eps,\Psi)$-solution  to the problem
\[
\tag{P} \text{minimize} \quad f(x) \quad \text{over} \quad x \in X.
\]
\end{definition}

Notice that, in view of Remark \ref{RemSolLinks}, the definition of $\Psi$-minimizing nets does not change if $(\eps,\Psi)$-solutions are replaced by $(\eps+,\Psi)$-solutions.

\begin{definition}
\label{DefWellPosed}
A function $f \colon X \to \P(Z,\Psi)$ is called $(0,\Psi)$-well-posed and $(0+,\Psi)$-well-posed, respectively, if a $(0,\Psi)$-solution and a $(0+,\Psi)$-solution, respectively, exists and every $\Psi$-minimizing net has a Hausdorff convergent subnet converging to a $(0,\Psi)$-solution and a $(0+,\Psi)$-solution, respectively of the problem (P).
\end{definition}

This definition is fundamentally different from the recent definition of well-posedness for set optimization problems as given in \cite{LiTeoZhang09NA, GutierrezEtAl12NA} which is closest in spirit to Definition \ref{DefWellPosed}. In fact, the concept due to Zhang, Li and Teo only deals with minimal elements with respect to the order relation $\lel_C$ on a normed space $Z$ defined through $A \lel_C B$ iff $B \subseteq A + C$ with a convex cone $C \subseteq Z$. On the other hand, in the same setting \cite{LongPeng13JOTA} considered well-posedness with respect to the order $\leu_C$ defined by $A \leu_C B$ iff $A \subseteq B - C$; no reasons are given why one or the other concept should be preferred. The definitions of `global well-posedness' in \cite[Definition~3.5-3.7]{MiglierinaMolhoRocca05JOTA} (see also \cite{LongPeng13JOTA}, for example) actually come closer to our definition (in particular, M-well-posedness) since they take the whole set of minimal points (of a vector-valued function) into account. 

The new features of Definition \ref{DefWellPosed} are the following: (i) it refers to sets of approximate solutions instead of single points and to the infimum in a complete lattice of sets which is completely absent in all previous concepts;  (ii) since the function $f$ maps into a complete lattice of sets, there is no ambiguity which ``set relation" has to be used, and it is clear that corresponding concepts for maximization have to be based on the sup-extension lattice (see Remark \ref{RemSupExt}); (iii) it covers many more set order relations and does not depend on the assumption that the ordering cone has a non-empty interior (as in many references, e.g.,  \cite[Definition 7-9]{CrespiKuroiwaRocca15AOR}).

The following example, taken from \cite[Remark 5.1]{MiglierinaMolhoRocca05JOTA} shows that Definition \ref{DefWellPosed} even recovers the missing ``scalarized" well-posedness in some cases.

\begin{example} Consider the set $S := \cb{x \in \R^2 \mid x_1 + x_2 \geq 0}$, the cone $C = \R^2_+$ with $\Psi = B^+ = \cb{w \in C^+ = \R^2_+ \mid w_1 + w_2 = 2}$ and define a function $f \colon \R^2 \to \P(\R^2, B^+)$ by
\[
f(x) = \left\{ 
	\begin{array}{ccc}
	x + \R^2_+ & \text{if} & x \in S \\
	\emptyset & \text{if} & x \not\in S
	\end{array}
	\right.
\]
Then, $\Min(f, 0, \Psi) = \cb{x \in \R^2 \mid x_1 + x_2 = 0}$, and this set also is a (full) lattice-solution of (P). Moreover, $\inf f[X] = \cb{z \in \R^2 \mid z_1 + z_2 \geq 0}$. Define $\psi_w(z) = w^Tz$ and hence $\psi^\triup_w(D) = \inf_{z \in D} w^Tz$. Then, $\psi^\triup_w(\inf f[X]) = -\infty$ if $w \neq (1,1)^T$ and $\psi^\triup_w(\inf f[X]) = 0$ if $w = (1,1)^T$. Hence, the set $M_\eps = \cb{x \in S \mid x_1 + x_2 \leq \eps}$ is the largest $(f, \eps, \Psi)$-minimizer for $\eps \geq 0$. The net $\cb{M_\eps}_{\eps > 0}$ is $\Psi$-minimizing and Hausdorff converges to $\Min(f, 0, \Psi)$. The same applies to $\Psi$-minimizing nets whose elements are subsets of $M_\eps$ for $\eps > 0$. Hence, $f$ is well-posed in the sense of Definition \ref{DefWellPosed} (using a set of linear $\psi$'s) whereas, according to \cite[Remark 5.1]{MiglierinaMolhoRocca05JOTA}, neither one of its linear scalarizations is (Tychonov) well-posed. 
\end{example}

It remains to be checked under what conditions the global well-posedness concepts for vector-valued or set-valued functions from \cite{MiglierinaMolhoRocca05JOTA} or \cite{LongPeng13JOTA} transfer into our concepts.

\begin{proposition}
\label{prop:WP}
Let $f \colon X \to \P(Z,\Psi)$ be given. If $f$ is $(0,\Psi)$-well-posed, then $\Min(f,0,\Psi)$ is non-empty and a $(0, \Psi)$-solution for (P). Moreover, $\Min(f,\eps,\Psi)$ Hausdorff converges to $\Min(f,0,\Psi)$ for $\eps \to 0$.
\end{proposition}

{\sc Proof.} Indeed, the assumptions imply the existence of a nonempty $(0,\Psi)$-solution for (P), hence $\Min(f,0,\Psi)$ is non-empty and also a $(0,\Psi)$-solution for (P). 

Since $\Min(f,\eps,\Psi) \neq \emptyset$ for $\eps > 0$ due to Proposition \ref{PropApproxSolExist} and $\Min(f,\eps_1,\Psi)\subseteq\Min(f,\eps_2,\Psi)$ for $0 \leq \eps_1 \leq \eps_2$, we only need to prove that for all neigborhoods $U$ of $0 \in X$ there exists a $\eps_U > 0$ such that $\Min(f,\eps,\Psi) \subseteq \Min(f,0,\Psi) + U$ is true for all $\eps < \eps_U$.

As $\cb{\Min(f,\eps,\Psi)}_{\eps\downarrow 0}$ is a $\Psi$-minimizing net and monotone with respect to inclusion, well-posedness implies the existence of a set $N\subseteq \Min(f,0,\Psi)$ such that for all $U$ it holds $\Min(f,\eps,\Psi)\subseteq N+U$ eventually. But this implies $\Min(f,\eps,\Psi)\subseteq \Min(f,0,\Psi)+U$ and thus Hausdorff convergence of $\Min(f,\eps,\Psi)$ to $\Min(f,0,\Psi)$ since $\Min(f,0,\Psi)\subseteq \Min(f,\eps,\Psi)+U$ is always true.
\pend

\medskip The above result remains true if $\Min(f,\eps,\Psi)$ is replaced by $\Min(f,\eps+,\Psi)$.

Within the framework of Corollary \ref{CorVectWeakEff}, $\WEff_{\eps e}(F,S)$ Hausdorff converges to $\WEff(F,S)\neq\emptyset$. Thus, Definition \ref{DefWellPosed} can be seen as a generalization of well-posedness concepts in vector optimization via well-posedness for the set-valued extension of a vector-valued function.

\begin{theorem}\label{thm:wp}
A function $f \colon X \to \P(Z,\Psi)$ is $(0,\Psi)$-well-posed if, and only if, the following three conditions are satisfied:

(1) $\Min(f,0,\Psi)\in \Sol(f,0,\Psi)$ is nonempty;

(2) Every net in $\Sol(f,0,\Psi)$ has a Hausdorff convergent subnet with a limit in $\Sol(f,0,\Psi)$;

(3) For all neigborhoods $U$ of $0 \in X$ there exists an $\eps>0$ such that $M_\eps \in \Sol(f, \eps, \Psi)$ implies the existence of $N \in\Sol(f,{\f 0,} \Psi)$ such that $M_\eps \subseteq N+U$ and $N \subseteq M_\eps+U$.
\end{theorem}

{\sc Proof.} First, assume that $f$ is $(0,\Psi)$-well-posed. Then a $(0,\Psi)$ solution exists and hence especially $\Min(f,0,\Psi)\in\Sol(f,0,\Psi)$ is satisfied. By assumption every net in $\Sol(f,0,\Psi)$ has a Hausdorff convergent subnet with limit in $\Sol(f,0,\Psi)$.
Finally, assume that for some neigborhoods $U$ of $0 \in X$ and all $n\in \N$ there exists a $M_n\in\Sol(f,\frac{1}{n},\Psi)$ such that for all $N\in\Sol(f,0,\Psi)$ it holds either $M_n\not\subseteq N+U$ or $N\not\subseteq M_n+U$. Especially, $\cb{M_n}_{n\in\N}$ is a $\Psi$-minimizing net and thus possesses a Hausdorff convergent subnet with limit $M\in\Sol(f,0,\Psi)$, contradicting the assumption.

On the other hand, if $\cb{M_i}_{i\in I}$ is a $\Psi$ minimizing net and $U$ a neighborhood of $0 \in X$, then eventually $M_i\in \Sol(f,\varepsilon,\Psi)$ is true for all $\varepsilon>0$ and it exists $N_i\in\Sol(f,0,\Psi)$ such that $M_i\subseteq N_i+U$ and $N_i\subseteq M_i+U$. But as  $\cb{N_i}_{i\in I}$ has a Hausdorff convergent subnet with limit $N\in\Sol(f,0,\Psi)$, so the same is satisfied for $\cb{M_i}_{i\in I}$.
\pend

\medskip Again, the result remains true if $\Sol(f, \eps, \Psi)$ is replaced by $\Sol(f, \eps+, \Psi)$ in (3) of the theorem.
Theorem \ref{thm:wp} generalizes a well-known characterization of generalized Tykhonov well-posedness for real-valued functions, compare \cite[Proposition 10.1.7]{Lucchetti06Book}.

\end{document}